\documentclass[preprint,10pt]{elsarticle}

\usepackage{amsmath,amssymb,amsfonts}
\usepackage{algorithmic}
\usepackage{graphicx}
\usepackage{textcomp}
\usepackage{color}
\usepackage{yhmath}                          % you prefer.
\usepackage{multicol}
\graphicspath{{figure/}}
\UseRawInputEncoding                               % you prefer.
\makeatletter

\newcommand{\Rmnum}[1]{\expandafter\@slowromancap\romannumeral #1@}
\makeatother
\usepackage{geometry}
\usepackage{picins}

\journal{Nuclear Physics B}

\begin{document}

\begin{frontmatter}


\title{Event-triggered adaptive consensus of heterogeneous multi-agent system under communication and actuator faults}

\author[Paestum]{Leyi Zheng}\ead{ly.zheng@siat.ac.cn}    % Add the
\author[Paestum]{Yimin Zhou}\ead{ym.zhou@siat.ac.cn}   

\address[Paestum]{Shenzhen Institute of Advanced Technology, Chinese Academy of Sciences, Shenzhen 518055, China and the University of Chinese Academy of Sciences, Beijing, China}  % Please supply     % here.

\begin{abstract}
In this paper, a heterogeneous leader-followers multiagent system is studied under simultaneous time-varying communication faults and actuator faults. First, the state of the leader is modelled as the closed-loop reference model where the states of the direct-connected followers are fed to the leader to improve the leader-followers tracking capability. An event-triggered communication mechanism is then designed for the agent information sharing among its neighbors so as to reduce the communication burden. Considering the time-varying communication link failure, a new distributed event-triggered observer is designed for each follower to estimate the whole system states so as to reduce the state error, whereas an adaptive distributed event-triggered estimator is further designed for the nondirect connected followers to estimate the coefficient matrix of the leader system. Further, an estimator is designed for the actuator fault estimation to reduce their impact on the system consistency. Hence, an adaptive event-triggered control strategy is proposed to ensure the consistency of the leader-follower system under the time-varying communication link faults and actuator faults. It is also shown that Zeno behavior is excluded for each agent and the effectiveness of the proposed adaptive event-triggered control strategy is verified on the heterogeneous multi-agent system.
\end{abstract}


\end{frontmatter}


\section{Introduction}
Multi-agent distributed cooperative control has received extensive attention in the past decades with its wide applications in smart grid, multi-robot system, intelligent transportation and other fields \cite{ref1}\cite{ref2}\cite{ref3}. In practice, multiple agents with limited perception capability should cooperate with each other to complete the task with specific formation \cite{ref4}. An interesting problem in multi-agent cooperative control is consistency \cite{ref5}, such as trajectory tracking \cite{ref6}, formation control\cite{ref7}, distributed estimation \cite{ref8}, motion planning \cite{ref9}. As the development of multi-agent systems, the system security has received more and more attention \cite{ref10}, \cite{ref11}. 

If the communication network is attacked, the communication links among the agents would be failed resulting in the inconsistency even failure of the task implementation \cite{ref12}. At present, many studies have analyzed the resilience of multi-agent systems against system faults and malicious attacks. In \cite{ref13}, a discrete control strategy is designed to solve the consistency problem of linear multi-agent systems under random link failures and communication noise. A random malicious node state estimation problem is discussed, where an orthogonal complement method is proposed to distinguish between malicious nodes and faulty agents \cite{ref14}. Then, the adverse effects on the attacked distributed multiple agents are analyzed to illustrate how the attack of a single node would spread to the entire network. A control method to mitigate the sensor attacks and actuator attacks is proposed in \cite{ref15} via the estimates of the expected behavior of the normal nodes. A distributed $H_{\infty}$ controller has been proposed and a distributed adaptive compensator is designed to reduce the interference of the attacks so as to tackle the problem of sensor and actuator failures \cite{ref16}. Most existing distributed controllers would lose function during communication interruption caused by the attacks. To solve this problem, a new hybrid distributed control strategy is proposed with the latest cached data to enhance the robustness of the multi-agent system \cite{ref17}.

As the number of agents grows, the system would experience a concomitant escalation in communication demands, thereby increasing the vulnerability to the communication link failures. Since the event-triggered mechanism uses discrete communication, it can effectively reduce the communication burden. For sensor and actuator attacks in the homogeneous systems, a distributed compensator is designed to predict the system state and an adaptive dynamic event-triggered mechanism (ETM) is proposed so that an elastic control strategy is designed to achieve the multi-agents consistency \cite{ref18}. In order to eliminate the measurement error, a self-triggering security control strategy is developed with the combination of system state and information of neighbor agents at the previous triggering interval \cite{ref19}. An event-triggered control strategy based on the state threshold is proposed in time-varying multi-agent systems to solve the problem of spurious data injection and parameter uncertainty \cite{ref20}. For the problem of malicious nodes in multi-agent homogeneous systems, a malicious nodes detection method is proposed via the neighbor information \cite{ref21}, and a distributed event-triggered mechanism $\&$ an adaptive event-triggered controller are developed to reduce the adverse effects of the malicious nodes \cite{ref22}. Paper \cite{ref23} has studied the unified consistency problem of the disturbed second-order multi-agent system, where a new framework of the event-based sliding mode control method is proposed with the distributed event-based sliding mode controllers and event-triggered detection. 

Further, the event-triggered mechanism with adjustable parameters is designed to better eliminate the Zeno phenomenon and reduce the amount of information transmission for a class of nonlinear multi-agent systems \cite{ref24}. In \cite{ref25}, the consistency of a linear heterogeneous multi-agent system is studied where the system matrix of the reference model is only known by some agents. An event-triggered adaptive distributed observer is designed for each node to estimate the system state and system matrix. Besides, the distributed event-triggered observer can be used to estimate the leader state under parameter uncertainties \cite{ref26}. In the uncertain nonlinear multi-agent tracking, the existing methods can only guarantee the bounded trajectory error tracking. In \cite{ref27}, followers are divided into two categories to design the distributed event-triggered observers for two groups of agents respectively. Based on this, the adaptive controllers are designed for subsystems with asymptotic convergence of all tracking errors. Although the followers are divided into two categories in \cite{ref25} and \cite{ref27}, the communication weights between different agents are not considered. Moreover, the two types of the followers cannot send messages to each other.

It is worth noting that most leaders only pass messages to the followers, however, the consistency performance of the multi-agent system may be improved with interactions between the leader and the followers\cite{ref28}, \cite{ref29}. In \cite{ref30}, a closed-loop reference model (CRM) is proposed to improve the transient performance of the system. Whereas an output feedback adaptive fault-tolerant control strategy has been proposed to track the CRM under external disturbances and actuator faults \cite{ref31}. The $H_{\infty}$ state tracking model reference adaptive control for switched systems is investigated \cite{ref32} with the CRM introduction to improve the transient performance. Then, a system with parameter uncertainty, external disturbance, actuator failure and other problems is discussed \cite{ref33}, identifying that the CRM can significantly improve the consistency of the system compared with the open-loop reference model (ORM). However, the above research has seldom been extended to multi-agent consistency systems. In \cite{ref34}, a distributed adaptive control with a CRM is applied in multi-agent CRM with the designed distributed parameter estimators and controllers for each agent. The simulation results show that the CRM has better instantaneous performance in multi-agent systems. 

\begin{figure}[!ht]
\centering
\includegraphics[width=0.4\textwidth]{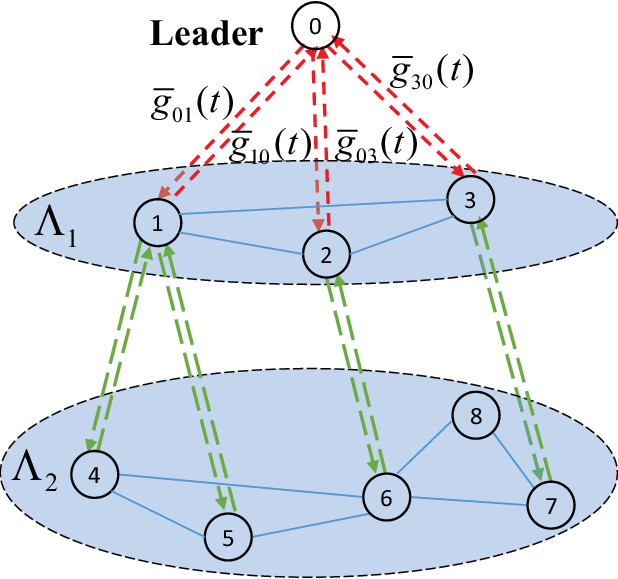}
\caption{The communication topology of the investigated heterogeneous multi-agent system.}
\label{fig_1}
%\vspace{-0.5em}
\end{figure}

Due to the complexity of the communication in multi-agent systems, communication faults and actuator faults may occur at the same time, where the consistency of the multi-agent systems is rarely discussed in the above situation. Hence, this paper propose an adaptive event-triggered observer and control strategy to achieve the consistency among the leader and the followers with the simultaneous communication faults and actuator faults. The followers are divided into two types: the first type is directly connected with the leader, and the second type has no direct connection with the leader, as depicted in Fig.\ref{fig_1}. The reference model receives information from the direct followers, and the leader and the followers in $\Lambda_1$ use a directed weighted communication method. Different from the setting in \cite{ref27}, the two types of the followers in this paper can communicate with each other via the weighted communication method triggered by event to reduce the communication burden.

The main contributions of the paper are summarized as follows.\\
\indent\setlength{\parindent}{1em}(1) An adaptive event-triggered consistency control strategy is developed for a heterogeneous multi-agent systems with time-varying communication faults and actuator faults to achieve the system consistency with uniform ultimate boundedness, while a CRM is established to receive the information from the first type followers. \\
\indent\setlength{\parindent}{1em} (2) The directed weighted communications are designed among the leader and the first type followers, and the first type followers $\&$ the second type followers. In order to reduce the communication burden, an adaptive event-triggered mechanism is designed, and the Zeno behavior can be excluded.\\
\indent\setlength{\parindent}{1em} (3) The event-triggered estimators are developed to estimate the leader state, follower states and actuator faults. The effectiveness of the proposed method is verified through simulation experiments.

The remainder of the paper is organized as follows. Section \Rmnum{2} describes the preliminary background and problem formulation. The main controller design procedure is explained in Section \Rmnum{3}. Section \Rmnum{4} is the simulation experiments and result analysis. Conclusion is provided in Section \Rmnum{5}.

\section{Preliminaries and Problem Formulation}
\subsection{Notion}
Let ${\mathbb R^{m \times n}}$ and ${\mathbb R ^n}$ denote the sets of $m \times n$ real matrices and $n$ dimensional real vectors; ${N_ + }$ represents a set of positive integers. Without exception, all the defined vectors are column vectors. ${I_n}$ is a n-dimensional identity matrix; ${1_n}$ and ${0_n}$ denote $n \times 1$ vectors with each entry as 1 and 0, respectively; $\otimes$ is the Kronecker product. For a real vector $x = [{x_1}, \cdots ,{x_N}] \in {\mathbb R^n}$, $||x||$ denotes its 2-norm $||x|{|_2}$. For a matrix $A \in {\mathbb R^{n \times n}}$, $||A||$ denotes the matrix norm $||A|{|_2}$. ${(.)^T}$ represents the transpose operator.
$\lambda (A)$ is the eigenvalue of $A$; ${\lambda _{\min }}(A)$ and ${\lambda _{\max }}(A)$ are the minimum and maximum of the eigenvalues. $diag({x_1}, \cdots ,{x_N})$ is a diagonal matrix with ${x_i}$, $i = 1, \cdots ,N$, on its diagonal position. $x \to 0$ denotes $x$ approaching zero. Let $vec([{x_1}, \cdots ,{x_N}]) = {[x_1^T, \cdots ,x_N^T]^T}$, $mat({[x_1^T, \cdots ,x_N^T]^T}) = [{x_1}, \cdots ,{x_N}]$, where ${x_i} \in {\mathbb R^n}$, $i = 1, \cdots ,N$. For brevity, $t$ is omitted in all variables hereafter for simplicity.

In this paper, the consistency between a leader and $N$ followers is studied. Assuming that the first $m$ followers are connected with the leader, described as ${\Lambda _1} = [1, \cdots m]$, while the rest followers are described as ${\Lambda _2} = [m + 1, \cdots N]$, without loss of generality.

\subsection{Graph theory}
The communications between the agents are modeled as a graph $G = (V,E,A)$, where $V = \{ 1, \cdots ,n\}$ is the node set, $E \in V \times V$ is the edge set of the pair agents, and $A = [{A_{11}}{\rm{  }}{A_{12}};{\rm{ }}{A_{21}}{\rm{  }}{A_{22}}] = [{a_{ij}}] \in {\mathbb R^{N \times N}}$ is the weighted adjacency matrix, where ${A_{11}}$ denotes the adjacency matrix among ${\Lambda _1}$, ${A_{22}}$ denotes the adjacency matrix among ${\Lambda _2}$, ${A_{12}}$ and ${A_{21}}$ denote the communications between ${\Lambda _1}$ and ${\Lambda _2}$. If node $j$ has an edge pointing to node $i$, node $j$ is the neighbor of node $i$ with ${a_{ij}}$ weight. The weights are defined as ${a_{ii}} = 0$, ${a_{ij}} > 0$ if $(j,i) \in E$, otherwise ${a_{ij}} = 0$. If ${a_{ij}} = {a_{ji}}$ for any pair of the neighboring agents, the associated communication graph is called an undirected graph, otherwise called a directed graph. Let $D = diag({d_i}, \cdots ,{d_n})$ denote the in-degree matrix, and ${d_i} = \sum\limits_{j = 1}^N {{a_{ij}}}$ is the weighted in-degree of node $i$. The Laplacian matrix of $G$ is defined as $L = D - A$. A sequence of edges $\{ ({i_1},{i_2}),({i_2},{i_3}), \cdots ,({i_{l - 1}},{i_l})\}$ is called a path from agent ${i_1}$ to agent ${i_l}$. The graph $G$ is a connected graph if for any two agents $i,j \in V$, there always exists a path from $i$ to $j$. Let ${N_L}$ be the collection of all direct followers to the leader (node 0). The symmetrical terms in a symmetric matrix are represented by $*$, that is, $[A \enspace * ;B\enspace C] = [A\enspace {B^T};B\enspace C]$.

\textbf{Assumption\,1:} Each follower has a path to the leader, that is, starting from any node $i\in V$, $i\in N_L$, there is a path to node $j\in V$, $j\notin N_L$.

\textbf{Assumption\,2:} The communications among the groups of the agents is undirected, that is, ${A_{11}}$ and ${A_{22}}$ are symmetric matrices. The leader $\&$ agents in ${\Lambda _1}$, and ${\Lambda_1}\& {\Lambda_2} $ both use directed communication strategy, that is, ${A_{12}}$ and ${A_{21}}$ are asymmetric, where ${A_{12}}$ represents the communication topology transmitting the information from ${\Lambda_2}$ to ${\Lambda_1}$, and ${A_{21}}$ is the communication topology from ${\Lambda _1}$ to ${\Lambda _2}$, as depicted in Fig.\ref{fig_1}.

\textbf{Remark\,1:} The $1^{st}$ type followers directly knows the leader's information and receives the communicated messages directly. The $2^{nd}$ type followers is unaware of the leader's information and relies on the $1^{st}$ type followers to receive the leader's messages. A estimator is designed to estimate the leader's coefficient matrix for the $2^{nd}$ type followers. The main purpose of the paper is to explore whether there is a significant difference in consistency between the two types of followers, when both actuator failures and communication failures occur in the system.

The communication link weight of each agent is defined as,
\begin{equation} \label{leader_1}
 \begin{array}{l}
{\bar a_{ij}} = {a_{ij}} + \delta _{ij}^a,
{\bar g_{0k}} = {g_{0k}} + \delta _{0k}^g,
{\bar g_{k0}} = {g_{k0}} + \delta _{k0}^g,
\end{array}
\end{equation}
where ${a_{ij}}$, ${g_{0k}}$ and ${g_{k0}}$ are the weights without communication fault; ${g_{0k}}$ is the communication weight of the leader to the direct neighbors in ${\Lambda _1}$; ${g_{i0}}$ is the communication weight of the followers in ${\Lambda _1}$ to the leader; $\delta _{ij}^a$, $\delta _{0k}^g$ and $\delta _{k0}^g$ are the unknown bounded communication disturbances, $i,j = 1, \cdots ,N$, $k = 1, \cdots ,m$.
\vspace{0.3\baselineskip}

\textbf{Remark\,2:} In this paper, the communication weights between the first type followers and leader may be different, that is, $g_{i0} \neq g_{0i}$, $i=1, \cdots, m$.

\begin{figure*}[ht]
\centering
\includegraphics[width=0.95\textwidth]{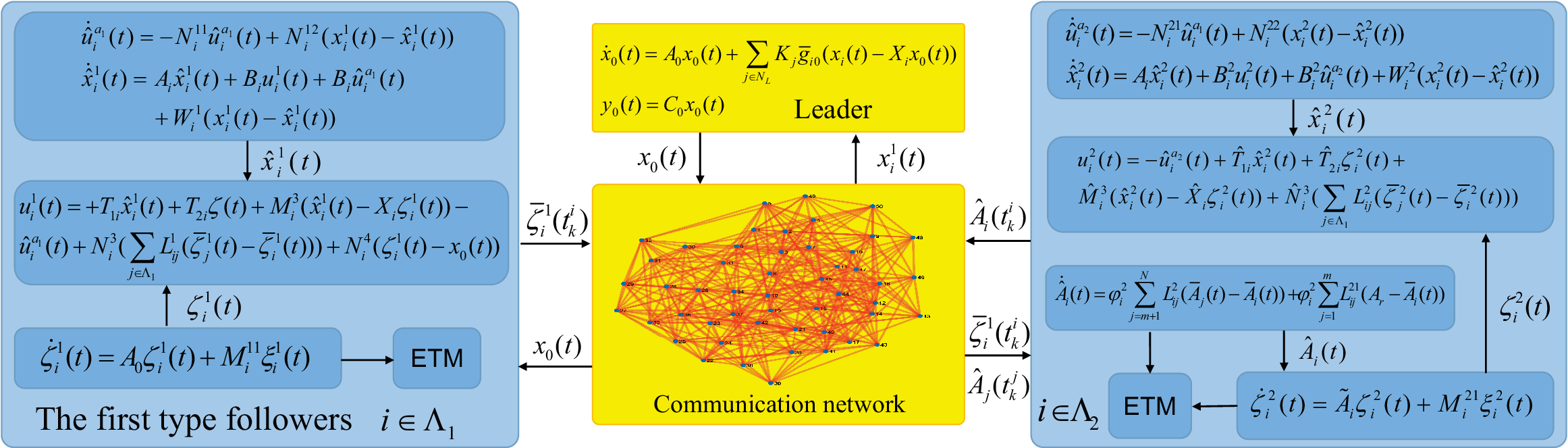}
\caption{The block diagram of the investigated multi-agent system.}
\label{control_1}
\end{figure*}

Let
$\bar A = [{\bar A_{11}}{\rm{  }} {\bar A_{12}};{\rm{ }}{\bar A_{21}}{\rm{  }} {\bar A_{22}}] = [{\bar a_{ij}}] \in {\mathbb R^{N \times N}}$, $L^{1} = {\bar A^{11}}$, $L^{2} = {\bar A^{22}}$,
$L^{12} = {\bar A^{12}}$, $L^{21} = {\bar A^{21}}$, $L_D^1 = {D^1} - {\bar A_{11}},{\rm{ }}L_D^2 = {D^2} - {\bar A_{22}}$,
${D^1} = diag([d_1^1, \cdots ,d_m^1]), {D^{12}} = diag([d_1^{12}, \cdots ,d_m^{12}]),$
${D^2} = diag([d_{m + 1}^2, \cdots ,d_N^2]),{D^{21}} = diag([d_{m + 1}^{21}, \cdots ,d_N^{21}]),$
${\bar G^1} = diag([{\bar g_{01}}, \cdots ,{\bar g_{0m}}]),{\bar G^2} = diag([{\bar g_{10}}, \cdots ,{\bar g_{m0}}]),$\\
where $d_i^1 = \sum\limits_{j = 1}^m {{{\bar a}_{ij}}}$, $d_i^{12} = \sum\limits_{j = m + 1}^N {{{\bar a}_{ij}}}$, $d_k^2 = \sum\limits_{j = m + 1}^N {{{\bar a}_{kj}}}$, $d_k^{21} = \sum\limits_{j = 1}^m {{{\bar a}_{kj}}}$, $i = 1, \cdots ,m$, $k = m + 1, \cdots ,N$.

\textbf{Definition\,1} If matrix $A$ is Hurwitz stable, there is a positive definite matrix $E$, so that $A^T E + EA$ is a negative definite matrix.
\vspace{0.3\baselineskip}

\textbf{Assumption\,3:} The communication weights ${\bar a_{ij}}>0,{\rm{ }} {\bar g_{0k}} > 0,{\rm{ }}{\bar g_{k0}} > 0$.

\subsection{Problem formulation}

Considering the following heterogeneous dynamic system with actuator fault,
\begin{equation}\label{leader_nof}
 \left\{
 \begin{array}{l}
{{\dot x}_i} = {A_i}{x_i} + {B_i}{u_i}\\
{y_i} = {C_i}{x_i}
\end{array}
\right.
\end{equation}
where ${x_i} \in {\mathbb R^{{n_i}}}$, ${u_i} \in {\mathbb R^{{m_i}}}$ and ${y_i} \in {\mathbb R^q}$ are the system state, control and system output, respectively. ${A_i}$, ${B_i}$ and ${C_i}$ are the real matrices with appropriate dimensions, $i = 1, \cdots ,N$. Here, the actuator fault is expressed as,
\begin{equation} \label{leader_2}
\begin{array}{l}
u_i = {u_i^b} + u_i^a,
\end{array}
\end{equation}
where $u_i^b$ is the normal control, $u_i^a$ is the unknown actuator fault of agent $i$.

\textbf{Assumption\,4:} The actuator fault $u_i^{{a}}$ and its derivative are unknown and bounded.

The reference model of the leader without feedback from the directed followers is described as,\\
\begin{equation}\label{leader_nof2}
\left\{
\begin{array}{l}
{{\dot {\hat {x} }}_0} = {A_0}{x_0}\\
{{\hat y}_0} = {C_0}{{\hat x}_0},
\end{array}
\right.
\end{equation}
where ${\hat x_0} \in {\mathbb R^{{p_1}}}$ is the leader state; ${\hat y_0} \in {\mathbb R^q}$ is the output of the reference model; ${A_0}$ and ${C_0}$ are the constant real coefficient matrices with appropriate dimensions.\\

The reference trajectory model with directed followers feedback is described as,
\begin{equation} \label{leader_nof3}
\left\{
\begin{array}{l}
{{\dot x}_0} = {A_0}{x_0} + \sum\limits_{j \in {N_L}} {{K_j}{{\bar g}_{i0}}({x_i} - {X_i}{x_0})} \\
{y_0} = {C_0}{x_0},
\end{array}
\right.
\end{equation}
where ${x_0} \in {\mathbb R^{{p_1}}}$ is the leader state; ${y_0} \in {\mathbb R^q}$ is the leader output; ${K_i}$ and ${X_i}$ are the gain real matrices with appropriate dimensions.

The aim of this paper is to design distributed observers-based controllers through event-triggered communication for System (\ref{leader_nof}), such that System (\ref{leader_nof}) under actuator failure and communication failure with reference closed-loop feedback System (\ref{leader_nof3}) can reach consistency with uniform ultimate boundedness(UUB), i.e., tracking errors $\{{y_i} - {y_0}\}$ is UUB.

Some assumptions are made during the solution of the above objective.

\textbf{Assumption\,5:} (1) $({A_0},{C_0})$ is observable; (2) $({A_i},{C_i})$ is controllable, $i = 1, \cdots ,N$.

\textbf{Assumption\,6:} The eigenvalues of ${A_0}$ are on the imaginary axis except zero point.

In this paper, the consistency between the leader and multiple followers in heterogeneous multi-agent systems is discussed. The pictorial diagram of the event-triggered observer and distributed control strategy for the two types of followers designed is shown in Fig.\ref{control_1}. The communication of the $1^{st}$ type followers with the leader is a weighted network, which can directly obtain the leader's state and affect the leader's state. The $2^{nd}$ type followers cannot communicate directly with the leader, but can communicate bidirectional with the $1^{st}$ type followers. This is different from the setting in \cite{ref27} where only the followers in ${\Lambda _1}$ affect the followers in ${\Lambda_2}$, however, $\Lambda_1$ and $\Lambda_2$ interact each other here, which would increase the difficulty in maintaining consistency of the multi-agent systems. The more the agents, the heavier the communication burden and higher difficulty of the system consistency to be achieved. Therefore, this paper adopts the event-triggered mechanism to transmit the system information so as to reduce the amount of the communication via the adaptive triggering degree of the agents. 

\section{Main results}
In order to reduce the impact of the actuator faults on the heterogeneous multi-agent systems, the actuator faults of the two types of the followers are estimated, and the observers are designed for each agent to estimate the system states.
\subsection*{A.1 Observer design for the  agents in ${\Lambda _1}$}
According to Assumption 4, the actuator faults estimator and the state observers of ${\Lambda _1}$ can be designed as,
\begin{equation} \label{leader_nof4}
\left\{
\begin{array}{l}
\dot {\hat u} _i^{{a_1}} =  - N_i^{11}\hat u_i^{{a_1}} + N_i^{12}(x_i^1 - \hat x_i^1)\\
\dot {\hat x}_i^1 = {A_i}\hat x_i^1 + {B_i}u_i^1 + {B_i}\hat u_i^{{a_1}} + W_i^1(x_i^1 - \hat x_i^1),
\end{array}
\right.
\end{equation}
where $\hat u_i^{{a_1}}$ is the estimate of the actuator fault $u_i^{{a_1}}$; $\hat x_i^1$ is the estimate of the system state $x_i^1$; $u_i^1$ is the controller of each agent in ${\Lambda _1}$; $N_i^{11}$, $N_i^{12}$ and ${W_i}$ are the matrices with appropriate dimensions, $i = 1, \cdots m$.

%According to system (\ref{leader_nof4}), it has,
Thus Eq.(\ref{leader_nof4}) can be rearranged as,
\begin{equation} \label{leader_nof5}
\left\{
\begin{array}{l}
{{\dot {\hat u}}^{{a_1}}} =  - {N^{11}}{{\hat u}^{{a_1}}} + {N^{12}}({x^1} - {{\hat x}^1})\\
{{\dot {\hat x}}^1} = {A^1}{{\hat x}^1} + {B^1}{u^1} + {B^1}{{\hat u}^{{a_1}}} 
+ {W^1}({x^1} - {{\hat x}^1}),
\end{array}
\right.
\end{equation}
where
${\hat u^{{a_1}}} = {[\hat u_1^{{a_1}T}, \cdots ,\hat u_m^{{a_1}T}]^T}$, ${x^1} = {[x_1^{1T}, \cdots ,x_m^{1T}]^T}$, ${\hat x^1} = {[\hat x_1^{1T}, \cdots ,\hat x_m^{1T}]^T}$, ${N^{11}} = diag\{ N_1^{11}, \cdots ,N_m^{11}\}$, ${N^{12}} = diag\{ N_1^{12}, \cdots ,N_m^{12}\} $, ${A^1} = diag\{ A_1^1, \cdots ,A_m^1\}$, ${B^1} = diag\{ B_1^1, \cdots ,B_m^1\}$ and ${W^1} = diag\{ W_1^1, \cdots ,W_m^1\}$.

Let
\begin{equation} \label{leader_nof5.1}
\begin{array}{l}
\tilde u_i^{{a_1}} = u_i^{{a_1}} - \hat u_i^{{a_1}}, \bar x_i^1 = x_i^1 - \hat x_i^1.
\end{array}
\end{equation}

%From system (\ref{leader_nof5}), it has,
Hence, Eq.(\ref{leader_nof5}) can be transformed as,
\begin{equation}
\begin{aligned}
{{\dot {\tilde {u}}}^{{a_1}}} =& {{\dot u}^{{a_1}}} + {N^{11}}{{\hat u}^{{a_1}}} - {N^{12}}({x^1} - {{\hat x}^1}),\\
{{\dot {\bar x} }^1}(t) =& ({A^1} - {W^1})({x^1} - {{\hat x}^1}) + {B^1}({u^{{a_1}}} - {{\hat u}^{{a_1}}}),
\end{aligned}
\end{equation}
where ${\tilde u^{{a_1}}} = {\left[ {\tilde u_1^{{a_1}T}, \cdots ,\tilde u_m^{{a_1}T}} \right]^T}$, ${\bar x^1} = {\left[ {\bar x_1^{1T}, \cdots ,\bar x_m^{1T}} \right]^T}$ and ${u^{{a_1}}} = {\left[ {u_1^{{a_1}T}, \cdots ,u_m^{{a_1}T}} \right]^T}$.

\textbf{Lemma\,1:} According to Assumption 4, appropriate gain matrices ${N^{11}}$, ${N^{12}}$ and ${W^1}$ can be set so that ${\tilde u^{{a_1}}}$ and ${\bar x^1}$ are UUB.

\textbf{Proof\,:} See Appendix A.

\subsection*{A.2 Observer design for the  agents in ${\Lambda _2}$}
According to Assumption 4, the actuator fault estimators and the state observer for each agent in ${\Lambda _2}$ can be designed as follows,
\begin{equation}
\begin{aligned}
\dot {\hat u}_i^{{a_2}} &= - N_i^{21}\hat u_i^{{a_2}} + N_i^{22}(x_i^2 - \hat x_i^2),\\
\dot {\hat x}_i^2 &= {A_i}\hat x_i^2 + B_i^2u_i^2 + B_i^2\hat u_i^{{a_2}} + W_i^2(x_i^2 - \hat x_i^2),
\end{aligned}
\end{equation}
where\\
${\hat u^{{a_2}}} = [\hat u_{m + 1}^{{a_2}}, \cdots ,\hat u_N^{{a_2}}]$, ${x^2} = [x_{m + 1}^2, \cdots ,x_N^2]$, ${\hat x^2} = [\hat x_{m + 1}^2, \cdots ,\hat x_N^2]$, ${N^{21}} = diag\{ N_{m + 1}^{21}, \cdots ,N_N^{21}\}$,\\ ${N^{22}} = diag\{ N_1^{22}, \cdots ,N_m^{22}\}$, ${A^2} = diag\{ A_{m + 1}^2, \cdots ,A_N^2\}$, ${B^2} = diag\{ B_{m + 1}^2, \cdots ,B_N^2\}$ and ${W^2} = diag\{ W_{m + 1}^2, \cdots ,W_N^2\}$.
Let
$${\tilde u^{{a_2}}} = {u^{{a_2}}} - {\hat u^{{a_2}}}, {\bar x^2} = {x^2} - {\hat x^2}.$$

Then, lemma 2 holds.

\textbf{Lemma\,2:} According to Assumption 4, appropriate gain matrices ${N^{21}}$, ${N^{22}}$ and ${W^2}$ can be set such that ${\tilde u^{{a_2}}}$ and ${\bar x^2}$ is UUB.

\textbf{Proof\,:} The proof is similar to Lemma 1.  $\hfill{} \Box$

\subsection*{B.1 Design of the distributed controller and trigger mechanism for the agents in ${\Lambda_1}$}

In order to reduce the communication burden, the communications among the  agents in ${\Lambda_1}$ are designed by the event-triggered mechanism.

The leader state observer is designed for the $1^{st}$ type followers,
\begin{equation}\label{type1_10.1}
\begin{aligned}
\dot \zeta _i^1=& {A_0}\zeta _i^1 + M_i^{11}\xi _i^1,\\
\xi _i^1=& M_i^{12}\sum\limits_{j \in \Lambda 1} {L_{ij}^1(\bar \zeta _j^1 - \bar \zeta _i^1)}  + {{\bar g}_{0i}}M_i^1({x_0} - \zeta _i^1)
+ M_i^{13}\sum\limits_{j = m + 1}^N {L_{ij}^{12}(\bar \zeta _j^2 - \bar \zeta _i^1)},
\end{aligned}
\end{equation}
where $\zeta _i^1$ represents the observation of the leader state by the ${i^{th}}$ follower; $M_i^{11}$, $M_i^{12}$, $M_i^1$ and $M_i^{13}$ are the constant matrices with compatible dimensions; $\zeta _j^2$ is the estimate of the leader state via the ${j^{th}}$ follower; $\bar \zeta _i^1$ and $\bar \zeta _j^2$ are the triggered states of the ${\Lambda_1}$ followers and ${\Lambda_2}$ followers based on the leader observation, respectively. For each agent $i$, $i \in {\Lambda _1}$, $\bar \zeta _i^1$ is defined as,
$$\bar \zeta _i^1 = {e^{{A_0}(t - t_k^i)}}\zeta _i^1(t_k^i),{\rm{  }}t \in [t_k^i,t_{k + 1}^i),$$
where $k \in {N_ + }$, $t_k^i$ is the ${k^{th}}$ triggering time instant when agent $i$ broadcasts the leader observation $\zeta _i^1$ to each of its neighbors.\\
Let
\begin{equation}
\begin{array}{l}
{\bar G^1} = diag([{\bar g_{01}}, \cdots ,{\bar g_{0m}}]),
{\bar G^2} = diag([{\bar g_{10}}, \cdots ,{\bar g_{m0}}]).
\end{array}
\end{equation}

\textbf{Remark\,3:} This paper considers the influence of the $2^{nd}$ type followers on the $1^{st}$ type followers. In such way, the state of the $2^{nd}$ type followers can be fed back to the leader indirectly through the $1^{st}$ type followers.

In order to introduce the event-triggered mechanism, the error functions are defined as,
$$e_{{v_i}}^1 = \bar \zeta _i^1 - \zeta _i^1,{\rm{ }}t \in [t_k^i,t_{k + 1}^i),{\rm{ }}i = 1, \cdots ,m.$$

The trigger functions are defined as,
\begin{small}
\begin{equation}
\begin{aligned}
{\dot \phi _{i1}} = & - {\beta _{i1}}{\phi _{i1}} - {\gamma _{i1}}(e_{{v_i}}^1{^T}{H_{i1}}e_{{v_i}}^1)
+ \frac{1}{{{\theta _{i1}}}}||{\psi _{i1}}|{|^2}\\
& + {\tau _{i1}}{e^{ - \frac{{{\delta _{i1}}}}{{{\sigma _{i1}} + t}}}} + \frac{{{{\hat \tau }_{i1}}}}{{{\hat\sigma _{i1}} + t}},
\end{aligned}
\end{equation}
\end{small}
where ${\psi _{i1}} = \sum\limits_{j \in {\Lambda _1}} {L_{ij}^1(\bar \zeta _j^1 - \bar \zeta _i^1)}$, ${\beta _{i1}} > 0$, ${\gamma _{i1}} > 0$, ${\tau _{i1}} > 0$, ${\delta _{i1}} > 0$, $\sigma _{i1}>0$, ${\hat \tau _{i1}}> 0$ and ${\hat \sigma _{i1}}> 0$ are the  parameters to be specified; ${H_{i1}}$ is a constant matrix with compatible dimension. Assuming that ${\phi _{i1}}(0)>0$, ${\theta _{i1}}$ is an adaptive parameter satisfying,
$${\theta _{i1}} = \frac{1}{{{d_i}}},{\rm{}}{d_i} = \sum\limits_{j = 1}^m {L_{ij}^1}.$$

It is noted that the first triggering time is $t_1^i = 0$. With the appropriate trigger mechanism, agent $i$ can determine the triggering time sequence $\{ l_k^i\} _{k = 2}^\infty$ by,
\begin{equation}\label{type1_10.2}
\begin{aligned}
l_{k + 1}^i =& \mathop {\inf }\limits_{r \ge t_k^i} \{ r:{\omega _{i1}}({\gamma _{i1}}(e_{{v_i}}^{1T}{H_{i1}}e_{{v_i}}^1)- {\theta _{i1}}||{\psi _{i1}}|{|^2} 
- {\tau _{i1}}{e^{ - \frac{{{\delta _{i1}}}}{{{\sigma _i} + t}}}} - \frac{{{{\hat \tau }_{i1}}}}{{{\hat\sigma _i} + t}})\ge {\phi _{i1}}\},
\end{aligned}
\end{equation}
where ${\omega _{i1}} > 0$ is the convergence function to be designed. Set ${\beta _{i1}}{\rm{ + }}\frac{1}{{{\omega _{i1}}}} > 0$, then ${\phi _{i1}} \ge {\phi _{i1}}(0){e^{ - \left( {{\beta _{i1}}{\rm{ + }}\frac{1}{{{\omega _{i1}}}}} \right)t}}$, $i \in {\Lambda _1}$.

\textbf{Remark\,4:} The designed triggering function can fulfill the requirements of the system consistency and effectively regulate the triggering intervals of the multi-agent system. For instances, $\frac{{{{\hat \tau }_{i1}}}}{{{\hat\sigma _i} + t}}$ is employed to govern the triggering intervals with relatively small $t$, and ${\tau _{i1}}{e^{ - \frac{{{\delta _{i1}}}}{{{\sigma _i} + t}}}}$ is utilized to regulate the triggering intervals with large $t$ instead. Here, $t$ is the time during which the multi-agent system is running.

The controllers of the $1^{st}$ type followers can be designed as,
\begin{equation}\label{type1_11}
\begin{aligned}
u_i^1 =&  - \hat u_i^{{a_1}} + {T_{1i}}\hat x_i^1 + {T_{2i}}\zeta  + M_i^3(\hat x_i^1 - {X_i}\zeta _i^1)
+ N_i^3(\sum\limits_{j \in {\Lambda _1}} {L_{ij}^1(\bar \zeta _j^1 - \bar \zeta _i^1)} ) + N_i^4(\zeta _i^1 - {x_0}),
\end{aligned}
\end{equation} 
where ${T_{1i}}$, ${T_{2i}}$, $M_i^3$, $N_i^3$ and $N_i^4$ are the constant matrices with appropriate dimensions, $i = 1, \cdots ,m$.

\textbf{Lemma\,3:} According to Assumption 4 and Assumption 5, an appropriate ${T_{1i}}$ is set so that the real part of the eigenvalues of matrix $({A_i} + {B_i}{T_{1i}})$ are less than zero. Let ${T_{2i}} = {Y_i} - {T_{1i}}{X_i}$, $({X_i},{Y_i})$ is the solution of the following equations,
\begin{equation}
\begin{aligned}
{X_i}{A_0} & = {A_i}{X_i} + {B_i}{Y_i} ,\\
{C_i}{X_i} & = {C_0},{\rm{   }}i = 1, \cdots ,m.
\end{aligned}
\end{equation}

\textbf{Lemma\,4:} According to Eq.(\ref{type1_11}) and Assumption 4, $N_i^3$ and $N_i^4$ are the matrices with appropriate dimensions, which can be set so that ${X_i}M_i^{11}M_i^{12} = {B_i}N_i^3$ and ${X_i}M_i^{11}M_i^1{\bar G^1} =  -{B_i}N_i^4$.

Let
\begin{equation}
\begin{array}{l}
{{\tilde x}^1} = {[{{\tilde x}_1}^{1T}, \cdots ,{{\tilde x}_m}^{1T}]^T},
{\rm{  }}{{\tilde x}^2} = {[{{\tilde x}_{m + 1}}^{2T}, \cdots ,{{\tilde x}_N}^{2T}]^T},
\end{array}
\end{equation}
\begin{equation}
\begin{array}{l}
{\varepsilon ^1} = {[\varepsilon _1^{1T}, \cdots ,\varepsilon _m^{1T}]^T},
{\rm{  }}e_v^1 = {[e_{{v_1}}^{1T}, \cdots ,e_{{v_m}}^{1T}]^T},
\end{array}
\end{equation}
\begin{equation}
\begin{array}{l}
e_v^2 = {[e_{{v_{m + 1}}}^{2T}, \cdots ,e_{{v_N}}^{2T}]^T},
\end{array}
\end{equation}

\noindent where $\tilde x_i^1 = \zeta _i^1 - {x_0}$, $\tilde x_j^2 = \zeta _j^2 - {x_0}$, $\varepsilon _i^1 = x_i^1 - {X_i}\zeta _i^1$, $i = 1, \cdots ,m$, $e_{{v_j}}^2 = \bar \zeta _j^2 - \zeta _j^2$, $\tilde x_j^2$ represents the state of agents $j$ in $\Lambda_2$, $j = m + 1, \cdots ,N$.

From Eqs.(\ref{leader_nof}), (\ref{leader_nof3}) and (\ref{type1_10.1}), it has,
\begin{equation}\label{type1_12}
\begin{aligned}
{{\dot {\tilde x}}^1} =& (({I_m} \otimes {A_0}) - {M^{11}}{M^1}({{\bar G}^1} \otimes {I_n}) - {M^{11}}{M^{12}}(L_D^1 \otimes {I_n}) \\
&- {M^{11}}{M^{13}}({D^{12}} \otimes {I_n})){{\tilde x}^1} - {M^{11}}{M^{13}}({D^{12}} \otimes {I_n})e_v^1 \\
&+ {M^{11}}{M^{13}}({A^{12}} \otimes {I_n})(e_v^2 + {{\tilde x}^2}) \\
&- \left( {{1_m} \otimes {K_1}\sum\limits_{j \in {N_L}} {{{\bar g}_{j0}}\left( {{x_j} - {X_j}{x_0}} \right)} } \right),\\
{{\dot \varepsilon }^1} =& ({A^1} + {B^1}{T_1}){{\hat x}^1} + ({B^1}{T_2} - X({I_m} \otimes {A_0})){\zeta ^1}\\
&+ {B^1}{M^3}({{\hat x}^1} - X{\zeta ^1}) + {W^1}({x^1}- {{\hat x}^1})- X{M^{11}}{M^{13}}\\
&{(\sum\limits_{j = m + 1}^N {L_{1j}^{12}{{(\zeta _j^2 - \zeta _1^1)}^T}} , \cdots ,\sum\limits_{j = m + 1}^N {L_{mj}^{12}{{(\zeta _j^2 - \zeta _m^1)}^T}} )^T},
\end{aligned}
\end{equation}
where
${M^{11}} = diag\{ M_1^{11}, \cdots ,M_m^{11}\}$, ${M^1} = diag{\rm{\{ }}M_1^1, \cdots ,M_m^1\}$, ${M^{12}} = diag{\rm{\{ }}M_1^{12}, \cdots ,M_m^{12}\}$, ${M^{13}} = diag{\rm{\{ }}M_1^{13}, \cdots ,M_m^{13}\}$, ${T_1} = diag({T_{11}}, \cdots ,{T_{1m}})$ and ${T_2} = diag({T_{21}}, \cdots ,{T_{2m}})$.

\textbf{Lemma\,5:} According to Assumptions 1, 2, 4, 5 and Definition 1, for the given parameter $a > 0$, the appropriate parameters ${c_1}$, ${c_2}$,  ${c_3}$ and matrix $E$, $Q$, $H$ can be set, and there is a positive definition matrix $P$, so that the inequality ${\Delta ^1} < 0$, and ${\Delta ^1} = [\Delta _{ij}^1] \in {\mathbb R^{3n \times 3n}}$ is established, where
$\Delta _{11}^1 = PE + {E^T}P + cP$, $\Delta _{21}^1 =  - {Q^T}P$, $\Delta _{22}^1 =  - {c_2}{\Upsilon _1}$, $\Delta _{31}^1 = \Delta _{32}^1 =  - L_D^{1T} \otimes {I_{\hat q}}$, $\Delta _{12}^1 = \Delta _{21}^{1T}$, $\Delta _{13}^1 = \Delta _{31}^{1T}$, $\Delta _{23}^1 = \Delta _{32}^{1T}$, $\Delta _{33}^1 =  - {c_3}{I_n}$, and ${\Upsilon _1}$ is a positive definition matrix.

\textbf{Proof\,:} It can be proved by diagonal dominance.  $\hfill{} \Box$

Considering the following Lyapunov function,
\begin{equation}\label{type1_12.1}
\begin{array}{l}
{V_{2i}} = V_{2i}^1 + V_{2i}^2,
{\rm{  }}V_{2i}^1 = \tilde x_i^{1T}{P_{1i}}{\tilde x^1} + {\phi _{i1}},
{\rm{  }}V_{2i}^2 = \varepsilon _i^{1T}{P_{2i}}\varepsilon _i^1,
\end{array}
\end{equation}

Let ${V_2} = \sum\limits_{j = 1}^m {{V_{2j}}}$, $V_2^1 = \sum\limits_{j = 1}^m {V_{2j}^1}$, $V_2^2 = \sum\limits_{j = 1}^m {V_{2j}^2}$. 

The time derivative of $V_2^1$ along the trajectory of System (\ref{leader_nof3}) is,
\begin{equation}
\begin{aligned}
\dot V_2^1 =& {{\tilde x}^{1T}}(E_1^T{P_1} + {P_1}{E_1}){{\tilde x}^1} - 2{{\tilde x}^{1T}}{P_1}{Q_1}e_v^1
+ 2{{\tilde x}^{1T}}{M^{11}}{M^{13}}({A^{12}} \otimes {I_n})e_v^2 - 2{{\tilde x}^{1T}}{K^1}{\Xi _1}\\
&+ 2{{\tilde x}^{1T}}{M^{11}}{M^{13}}({A^{12}} \otimes {I_n}){{\tilde x}^2} + \sum\limits_{j = 1}^m ({\hat\delta_{j1}} - {{\beta _{j1}}} {\phi _{j1}})
- \sum\limits_{j = 1}^m {({\gamma _{j1}}(e_{{v_j}}^1{{}^T}{H_{j1}}e_{{v_j}}^1) - \frac{1}{{{\theta _{j1}}}}||{\psi _{j1}}|{|^2})},
\end{aligned}
\end{equation}
where\\
${E_1} = (({I_m} \otimes {A_0}) - {M^{11}}{M^1}({\bar G^1} \otimes {I_n}) - {M^{11}}{M^{12}}(L_D^1 \otimes {I_n}) - {M^{11}}{M^{13}}({D^{12}} \otimes {I_n})),$\\
${Q_1} = {M^{11}}{M^{13}}({D^{12}} \otimes {I_n}),$ 
${\Xi _1} = {1_m} \otimes \sum\limits_{j \in {N_L}} {{{\bar g}_{j0}}({x_j} - {X_j}{x_0})}$, $\hat\delta_{j1}={\tau _{j1}}{e^{ - \frac{{{\delta _{j1}}}}{{{\sigma _{j1}} + t}}}} + \frac{{{{\hat \tau }_{j1}}}}{{{\hat\sigma _{j1}} + t}}$.

Let ${\theta _{\min }} = \min \{ {\theta _{11}}, \cdots ,{\theta _{m1}}\}$, ${\gamma _{\min }} = \min \{ {\gamma _{11}}, \cdots ,{\gamma _{m1}}\}$, ${H_1} = diag({H_{11}}, \cdots ,{H_{m1}})$.

Then
\begin{equation}
\begin{array}{l}
\sum\limits_{j = 1}^m {\frac{1}{{{\theta _{i1}}}}||{\psi _{i1}}|{|^2}}  \le \frac{1}{{{\theta _{\min }}}}\sum\limits_{j = 1}^m {||{\psi _{i1}}|{|^2}} 
 = \frac{1}{{{\theta _{\min }}}}\chi _1^T{[{I_n},{I_n}]^T}(L_D^{1T}L_D^1 \otimes {I_n})[{I_n},{I_n}]{\chi _1},
\end{array}
\end{equation}
where ${\chi _1} = {[{\tilde x^{1T}},e_v^{1T}]^T}$.

Then, it has,
\begin{equation}\label{type1_13}
\begin{array}{l}
{\tilde x^{1T}}(E_1^T{P_1} + {P_1}{E_1}){\tilde x^1} - 2{\tilde x^{1T}}{P_1}{Q_1}e_v^1 
- \sum\limits_{j = 1}^m {({\gamma _{i1}}e_{{v_i}}^1{{}^T}{H_{i1}}e_{{v_i}}^1)} \le \chi _1^T{\Gamma ^1}\chi _1^T,
\end{array}
\end{equation}
\begin{equation}\label{type1_14}
\begin{array}{l}
 - 2{{\tilde x}^{1T}}{K^1}{\Xi _1}
 =  - 2{{\tilde x}^{1T}}{K^1}({1_m} \otimes \sum\limits_{j \in {N_L}} {{g_{j0}}(x_j^1 - \hat x_j^1))} 
 - 2{{\tilde x}^{1T}}{K^1}({1_m} \otimes \sum\limits_{j \in {N_L}} {{g_{j0}}(\hat x_j^1 - {X_j}\zeta _j^1 + {X_j}\zeta _j^1 - {X_j}{x_0}))},
\end{array}
\end{equation}
where ${\Gamma ^1} = [\Gamma _{11}^1,\Gamma _{12}^1;*,\Gamma _{22}^1]$ with $\Gamma _{11}^1 = E_1^T{P_1} + {P_1}{E_1} + {c_1}{P_1}$, $\Gamma _{21}^1 =  - Q_1^T{P_1}$, $\Gamma _{22}^1 =  - {\gamma _{\min }}{H_1}$ and ${c_1} > 0$.

Set an appropriate matrix ${K^1}$, so that the following inequality holds,
\begin{equation}\label{type1_15}
\begin{array}{l}
- 2{\tilde x^{1T}}{K^1}({1_m} \otimes \sum\limits_{j \in {N_L}} {{g_{j0}}({X_j}\zeta _j^1 - {X_j}{x_0}))} < 0.
\end{array}
\end{equation}

By Young's inequality,
\begin{equation}\label{type1_16}
\begin{array}{l}
2{\tilde x^{1T}}{M^{11}}{M^{13}}({A^{12}} \otimes {I_n})e_v^2 
\le {\beta _1}{\tilde x^{1T}}{\tilde x^1} + \frac{1}{{{\beta _1}}}||{M^{11}}{M^{13}}({A^{12}} \otimes {I_n})e_v^2|{|^2},
\end{array}
\end{equation}
\begin{equation}\label{type1_17}
\begin{array}{l}
2{\tilde x^{1T}}{M^{11}}{M^{13}}({A^{12}} \otimes {I_n}){\tilde x^2}
\le {\beta _2}{\tilde x^{1T}}{\tilde x^1} + \frac{1}{{{\beta _2}}}||{M^{11}}{M^{13}}({A^{12}} \otimes {I_n}){\tilde x^2}|{|^2},
\end{array}
\end{equation}
\begin{equation}\label{type1_18}
\begin{array}{l}
- 2{\tilde x^{1T}}{K^1}({1_m} \otimes \sum\limits_{j \in {N_L}{\beta _3}} {{g_{j0}}(x_j^1 - \hat x_j^1)} )
\le {\beta _3}||{\tilde x^{1T}}{K^1}|{|^2} + \frac{1}{{{\beta _3}}}||{\bar g_1}{\hat K_1}{\bar x^1}|{|^2},
\end{array}
\end{equation}
\begin{equation} \label{type1_19}
\begin{array}{l}
- 2{\tilde x^{1T}}{K^1}({1_m} \otimes \sum\limits_{j \in {N_L}} {{g_{j0}}(\hat x_j^1 - {X_j}\zeta _j^1))} 
\le {\beta _4}||{\tilde x^{1T}}{K^1}|{|^2} + \frac{1}{{{\beta _4}}}||{\bar g_1}{\hat K_1}{\varepsilon ^1}|{|^2},
\end{array}
\end{equation}
where ${\hat K_1} = ({1_m} \otimes 1_m^T) \otimes {I_n}$; ${\bar g_1} = \max \{ {\bar g_{10}}, \cdots ,{\bar g_{m0}}\}$; ${\beta_1}$, ${\beta_2}$, ${\beta_3}$ and ${\beta_4}$ are the real values to be determined.
According to Eqs.(\ref{type1_13})-(\ref{type1_19}), Lemma 5 and Schur Complement Lemma, it has,
\begin{small}
\begin{equation} \label{type1_20}
\begin{aligned}
\dot V_2^1 \le &  - {{\hat \alpha }_1}{c_1}{{\tilde x}^{1T}}{P_1}{{\tilde x}^1} - {{\hat \alpha }_2}{c_1}{{\tilde x}^{1T}}{P_1}{{\tilde x}^1} + {\beta _1}{{\tilde x}^{1T}}{{\tilde x}^1} - \sum\limits_{j = 1}^m {{\beta _{j1}}} {\phi _{j1}}\\
& + {\beta _2}||{M^{11}}{M^{13}}({A^{12}} \otimes {I_n}){{\tilde x}^1}|{|^2}
+ {\beta _3}||{{\tilde x}^{1T}}{K^1}|{|^2} + \frac{1}{{{\beta _2}}}{{\tilde x}^{2T}}{{\tilde x}^2}\\
&+ \sum\limits_{j = 1}^m {\hat\delta_{j1}} + {\beta _4}||{{\tilde x}^{1T}}{K^1}||^2 + \frac{1}{{{\beta _4}}}||{{\bar g}_1}{{\hat K}_1}{\varepsilon ^1}|{|^2} + {\Omega _1}, \\
\end{aligned}
\end{equation}\end{small}
where ${\Omega _1} = \frac{1}{{{\beta _1}}}||{M^{11}}{M^{13}}({A^{12}} \otimes {I_n})e_v^2|{|^2} + \frac{1}{{{\beta _3}}}||{\bar g_1}{\hat K_1}{\bar x^1}|{|^2}$, ${\hat \alpha _2} = 1 - {\hat \alpha _1}$ and $0 < {\hat \alpha _1} < 1$.
According to Lemma 1, ${\Omega _1}$ is bounded.

Let
\begin{equation}
\begin{aligned}
{\Omega _2} =&  - {\hat \alpha _1}{c_1}{\tilde x^{1T}}{P_1}{\tilde x^1} + {\beta _1}{\tilde x^{1T}}{\tilde x^1} + {\beta _3}||{\tilde x^{1T}}{K^1}|{|^2} \\
& + {\beta _2}||{M^{11}}{M^{13}}({A^{12}} \otimes {I_n}){\tilde x^1}|{|^2} 
  + {\beta _4}||{\tilde x^{1T}}{K^1}||^2.
\end{aligned}
\end{equation}

From Eq.(\ref{type1_20}), it has,
\begin{equation}\label{type1_20.1}
\begin{aligned}
\dot V_2^1 \le &  - {\hat \alpha _2}{c_1}{\tilde x^{1T}}{P_1}{\tilde x^1} + {\Omega _1} + {\Omega _2} 
+ \frac{1}{{{\beta _4}}}||{\bar g_1}{\hat K_1}{\varepsilon ^1}|{|^2} \\
& + \frac{1}{{{\beta _2}}}{\tilde x^{2T}}{\tilde x^2} - \sum\limits_{j = 1}^m {{\beta _{j1}}} {\phi _{j1}}
+ \sum\limits_{j = 1}^m {\hat\delta_{j1}}.
\end{aligned}
\end{equation}

\textbf{Lemma\,6:} According to Assumptions 1, 4 and Definition 1, an appropriate ${M^3}$ can be found so that there is a positive definite matrix ${P_2}$,
\begin{small}$${\left( {{A^1} + {B^1}({T_1} + {M^3})} \right)^T}{P_2} + {P_2}\left( {{A^1} + {B^1}({T_1} + {M^3})} \right) \le  - {\tilde \alpha _1}{P_2},$$ \end{small}where ${\tilde \alpha _1} > 0$. According to Lemma 3, the time derivation of $V_2^2$ along the trajectory of System (\ref{type1_12}) is obtained as,
\begin{small}
\begin{equation}
\begin{aligned}
\dot V_2^2 =& {\varepsilon ^1}^T{\left( {{A^1} + {B^1}({T_1} + {M^3})} \right)^T}{P_2}{\varepsilon ^1} + 2{\varepsilon ^1}^T{P_2}{W^1}{{\bar x}^1} + {\varepsilon ^1}^T{P_2}\left( {{A^1} + {B^1}({T_1} + {M^3})} \right){\varepsilon ^1} - 2{\varepsilon ^1}^T{P_2}X{M^{11}}{M^{13}}\\
& {(\sum\limits_{j = m + 1}^N {L_{1j}^{12}{{(\bar \zeta _j^2 - \bar \zeta _1^1)}^T}} , \cdots , \sum\limits_{j = m + 1}^N {L_{1j}^{12}{{(\bar \zeta _j^2 - \bar \zeta _m^1)}^T}} )^T}.
\end{aligned}
\end{equation}
\end{small}

Let
\begin{equation}
\begin{aligned}
h = & - 2{\varepsilon ^1}^T{P_2}X{M^{11}}{M^{13}} (\sum\limits_{j = m + 1}^N {L_{1j}^{12}{{(\bar \zeta _j^2 - \zeta _j^2 + \zeta _j^2 - {x_0} + {x_0} - \zeta _1^1 + \zeta _1^1 - \bar \zeta _1^1)}^T}} , \cdots ,\\
& \sum\limits_{j = m + 1}^N {L_{mj}^{12}{{(\bar \zeta _j^2 - \zeta _j^2 + \zeta _j^2 - {x_0} + {x_0} - \zeta _m^1 + \zeta _m^1 - \bar \zeta _m^1)}^T}} {)^T}.
\end{aligned}
\end{equation}

By Young's inequality,
\begin{small}
\begin{equation}\label{type1_21}
\begin{aligned}
{h_1} &= - 2{\varepsilon ^1}^T{P_2}X{M^{11}}{M^{13}}{\left( {\sum\limits_{j = m + 1}^N {L_{1j}^{12}e_{{v_j}}^{2T}} , \cdots ,\sum\limits_{j = m + 1}^N {L_{mj}^{12}e_{{v_j}}^{2T}} } \right)^T}
& \le {\tau _1}{\varepsilon ^1}^T{\varepsilon ^1} + \frac{1}{{{\tau _1}}}||{P_2}X{M^{11}}{M^{13}}({A^{12}} \otimes {I_n})e_v^2|{|^2},
\end{aligned}
\end{equation}
\begin{equation}\label{type1_22}
\begin{aligned}
{h_2} &= - 2{\varepsilon ^1}^T{P_2}X{M^{11}}{M^{13}}{\left( {\sum\limits_{j = m + 1}^N {L_{1j}^{12}\bar x_j^{2T}} , \cdots ,\sum\limits_{j = m + 1}^N {L_{mj}^{12}\bar x_j^{2T}} } \right)^T}
& \le {\tau _2}{\varepsilon ^1}^T{\varepsilon ^1} + \frac{1}{{{\tau _2}}}||{P_2}X{M^{11}}{M^{13}}({A^{12}} \otimes {I_n}){{\tilde x}^2}|{|^2},
\end{aligned}
\end{equation}
\begin{equation}\label{type1_23}
\begin{aligned}
{h_3} &= 2{\varepsilon ^1}^T{P_2}X{M^{11}}{M^{13}}{\left( {\sum\limits_{j = m + 1}^N {L_{1j}^{12}\tilde x{{_1^{1T}}}} , \cdots ,\sum\limits_{j = m + 1}^N {L_{mj}^{12}\tilde x_m^{1T}} } \right)^T}
&\le {\tau _3}{\varepsilon ^1}^T{\varepsilon ^1} + \frac{1}{{{\tau _3}}}||{P_2}X{M^{11}}{M^{13}}({A^{12}} \otimes {I_n}){{\tilde x}^1}|{|^2},
\end{aligned}
\end{equation}
\begin{equation} \label{type1_24}
\begin{aligned}
{h_4} &= - 2{\varepsilon ^1}^T{P_2}X{M^{11}}{M^{13}}{\left( {\sum\limits_{j = m + 1}^N {L_{1j}^{12}e_{{v_j}}^{1T}} , \cdots ,\sum\limits_{j = m + 1}^N {L_{mj}^{12}e_{{v_j}}^{1T}} } \right)^T}
& \le {\tau _4}{\varepsilon ^1}^T{\varepsilon ^1} + \frac{1}{{{\tau _4}}}||{P_2}X{M^{11}}{M^{13}}({A^{12}} \otimes {I_n})e_v^1|{|^2},
\end{aligned}
\end{equation}
\begin{equation} \label{type1_25}
\begin{array}{l}
2{\varepsilon ^1}^T{P_2}{W^1}{\bar x^1} \le {\tau _5}{\varepsilon ^1}^T{\varepsilon ^1} + \frac{1}{{{\tau _5}}}||{P_2}{W^1}{\bar x^1}|{|^2},
\end{array}
\end{equation}
\end{small}where ${\tau _1},{\tau _2},{\tau _3},{\tau _4},{\tau _5} > 0$, $\bar x_i^1 = \zeta _i^1 - {x_0}$, $\bar x_j^2 = \zeta _j^2 - {x_0}$, $i = 1, \cdots m$ and $j = m + 1, \cdots N$.

According to Lemma 6, Eqs.(\ref{type1_21})-(\ref{type1_25}), it has,
\begin{small}
\begin{equation} \label{type1_26}
\begin{aligned}
\dot V_2^2 =& {\tau _6}{\varepsilon ^1}^T{\varepsilon ^1} + \frac{1}{{{\tau _2}}}||{P_2}X{M^{11}}{M^{13}}({A^{12}} \otimes {I_n}){{\tilde x}^2}|{|^2} - {{\tilde \alpha }_1}{\varepsilon ^1}^T{P_2}{\varepsilon ^1} + \hat M + \frac{1}{{{\tau _3}}}||{P_2}X{M^{11}}{M^{13}}({A^{12}} \otimes {I_n}){{\tilde x}^1}|{|^2},
\end{aligned}
\end{equation}
\end{small}where $\hat M = \frac{1}{{{\tau _1}}}||{P_2}X{M^{11}}{M^{13}}({A^{12}} \otimes {I_n})e_v^2|{|^2} + \frac{1}{{{\tau _4}}}||{P_2}X{M^{11}}{M^{13}}({A^{12}} \otimes {I_n})e_v^1|{|^2} + \frac{1}{{{\tau _5}}}||{P_2}{W^1}{\bar x^1}|{|^2}$, ${\tau _6} = {\tau _1} + {\tau _2} + {\tau _3} + {\tau _4} + {\tau _5}$.

According to Lemma 1 and Eq.(\ref{type1_10.2}), it is known that $\hat M$ is bounded.

From Eqs.(\ref{type1_12.1}), (\ref{type1_20.1}) and  (\ref{type1_26}), it has,
\begin{small}
\begin{equation} \label{type1_27}
\begin{aligned}
{{\dot V}_2} =& \dot V_2^1 + \dot V_2^2
\le  - {{\hat \alpha }_2}{c_1}{{\tilde x}^{1T}}{P_1}{{\tilde x}^1} + {\Omega _3} + \frac{1}{{{\beta _2}}}{{\tilde x}^{2T}}{{\tilde x}^2} + \hat M + {\Omega _1}\\
&+ \frac{1}{{{\tau _2}}}||{P_2}X{M^{11}}{M^{13}}({A^{12}} \otimes {I_n}){{\tilde x}^2}|{|^2}
- {{\tilde \alpha }_1}{\varepsilon ^1}^T{P_2}{\varepsilon ^1} \\
&+ {\tau _6}{\varepsilon ^1}^T{\varepsilon ^1}
+ \frac{1}{{{\beta _4}}}||{{\bar g}_1}{{\hat K}_1}{\varepsilon ^1}|{|^2} 
- \sum\limits_{j = 1}^m {{\beta _{j1}}} {\phi _{j1}} + \sum\limits_{j = 1}^m {\hat\delta_{j1}},
\end{aligned}
\end{equation}
\end{small}where ${\Omega _3} = {\Omega _2} + \frac{1}{{{\tau _3}}}||{P_2}X{M^{11}}{M^{13}}({A^{12}} \otimes {I_n}){\tilde x^1}|{|^2}$.

\subsection*{B.2 Design of the distributed controller and trigger mechanism for the agents in ${\Lambda_2}$}

Here, the distributed event-triggered mechanism and event-triggered controllers are proposed for the $2^{nd}$ type followers. It is aware that the agents in $\Lambda_2$ have no knowledge of the leader state, so the leader state is estimated by each agent in $\Lambda_2$ as follows,
\begin{equation} \label{type2_1}
\begin{array}{l}
{\dot {\hat A}_i} = \varphi _i^2\sum\limits_{j = m + 1}^N {L_{ij}^2({{\bar A}_j} - {{\bar A}_i}) + } \varphi _i^2\sum\limits_{j = 1}^m {L_{ij}^{21}({A_r} - {{\bar A}_i})},
\end{array}
\end{equation}
where ${A_r} = vec({A_0})$; ${\hat A_i}$ is the estimation of ${A_r}$ by agent $i$; $\varphi _i^2$ is a constant real value to be determined, $i = m + 1, \cdots ,N$. For each agent $i$, $i = m + 1, \cdots ,N$, ${\bar A_i}$ is defined as,
$${\bar A_i} = {\hat A_i}(\hat t_k^i),{\rm{  }}t \in [\hat t_k^i,\hat t_{k + 1}^i),$$
where $k \in {N_ + }$, $\hat t_k^i$ is the $k^{th}$ triggering time instant when agent $i$ broadcasts the ${\hat A_i}$ via agent $i$.

Define the following triggering function,
\begin{equation} \label{type2_2}
\begin{aligned}
{\dot \phi _{i2}} =& - {\beta _{i2}}{\phi _{i2}} - {\gamma _{i2}}({\bar e_i}{^T}{\bar e_i}) + \frac{1}{{{\theta _{i2}}}}||{\psi _{i2}}|{|^2}
 + {\tau _{i2}}{e^{ - \frac{{{\delta _{i2}}}}{{{\sigma _{i2}} + t}}}} + {{\hat \tau }_{i2}},
\end{aligned}
\end{equation}
where ${\bar e_i} = {\hat A_i} - {\bar A_i}$, ${\psi _{i2}} = \sum\limits_{j \in {\Lambda _2}} {L_{ij}^2({{\bar A}_j} - {{\bar A}_i})}  + \sum\limits_{j = 1}^m {L_{ij}^{21}({A_r} - {{\bar A}_i})}$, ${\beta _{i2}} > 0$, ${\gamma _{i2}} > 0$, ${\theta _{i2}} > 0$, ${\tau _{i2}} > 0$, ${\delta _{i2}} > 0$, ${\sigma _{i2}}>0$ and ${\hat\tau _{i2}}>0$ are the parameters to be specified, $i = m + 1, \cdots ,N$.

The first triggering time is $\hat t_1^i = 0$, and agent $i$ determines the triggering time sequence $\{ \hat l_k^i\} _{k = 2}^\infty $ by,
\begin{equation} \label{type2_3}
\begin{array}{l}
\hat l_{k + 1}^i = 
\mathop {\inf }\limits_{r \ge t_k^i} \{ \hat r:{\hat \omega _i}({\gamma _{i2}}||{\bar e_i}||^2 - \frac{1}{{{\theta _{i2}}}}||{\psi _{i2}}|{|^2} - {\tau _{i2}}{e^{ - \frac{{{\delta _{i2}}}}{{{\sigma _{i2}} + t}}}} - {\hat\tau _{i2}}) \ge {\phi _{i2}}\},
\end{array}
\end{equation}
where ${\hat \omega _i} > 0$ is the convergence function to be determined.

Let
\begin{equation}
\begin{array}{l}
\dot {\hat A} = {\varphi ^2}({H_2} \otimes {I_{{n^2}}})({\bar A_r} - \hat A + \hat A - \bar A) 
= {\varphi ^2}({H_2} \otimes {I_{{n^2}}})(\mathord{\buildrel{\lower3pt\hbox{$\scriptscriptstyle\smile$}}
\over A}  + \bar e),
\end{array}
\end{equation}
where ${\varphi ^2} = diag(\varphi _{m + 1}^2, \cdots ,\varphi _n^2) \otimes {I_{{n^2}}}$, ${H_2} = L_D^2 + {D^{21}}$, ${\bar A_r} = diag({A_r}, \cdots ,{A_r})$, $\hat A = diag({\hat A_{m + 1}}, \cdots ,{\hat A_N})$, $\bar A = diag({\bar A_{m + 1}}, \cdots ,{\bar A_N})$, $\mathord{\buildrel{\lower3pt\hbox{$\scriptscriptstyle\smile$}}
\over A}  = {\bar A_r} - \hat A$ and $\bar e = {[\bar e_{m + 1}^T, \cdots ,\bar e_N^T]^T} = \hat A - \bar A$.
Then
\begin{equation} \label{type2_4}
\begin{array}{l}
\dot {\mathord{\buildrel{\lower3pt\hbox{$\scriptscriptstyle\smile$}}
\over A}}  =  - {\varphi ^2}({H_2} \otimes {I_{{n^2}}})(\mathord{\buildrel{\lower3pt\hbox{$\scriptscriptstyle\smile$}}
\over A}  + \bar e). 
\end{array}
\end{equation}

\textbf{Lemma\,7:} Considering System (\ref{type2_1}), ${\mathord{\buildrel{\lower3pt\hbox{$\scriptscriptstyle\smile$}}
\over A} _i} = ({A_r} - {\hat A_i})$ reaches the consistency with UUB under triggering function (\ref{type2_2}) and triggering mechanism (\ref{type2_3}), $i = m + 1, \cdots ,N$.

\textbf{Proof\,:} See Appendix B.

A leader state observer can be designed for the $2^{nd}$ type followers,
\begin{equation} \label{type2_19}
\begin{array}{l}
\dot \zeta _i^2 = {{\tilde A}_i}\zeta _i^2 + M_i^{21}\xi _i^2,\\
\xi _i^2 = M_i^{22}\sum\limits_{j \in {\Lambda _1}} {L_{ij}^2(\bar \zeta _j^2 - \bar \zeta _i^2)}
+ M_i^{23}\sum\limits_{j = 1}^m {L_{ij}^{21}(\bar \zeta _j^1 - \bar \zeta _i^2)},
\end{array}
\end{equation}
where $\zeta _i^2$ represents the estimate of the leader state by the $i^{th}$ follower; $M_i^{21}$, $M_i^{22}$ and $M_i^{23}$ are the constant matrices with compatible dimensions; ${\tilde A_i} = mat({\hat A_i})$, $i \in {\Lambda _2}$; $\bar \zeta _i^2$ is defined as, $\bar \zeta _i^2 = \zeta _i^2(\tilde t_k^i),{\rm{  }}t \in [\tilde t_k^i,\tilde t_{k + 1}^i)$, $\tilde t_k^i$ is the $k^{th}$ triggering time instant when agent $i$ broadcasts the leader observation status $\zeta _i^2$ to each of its neighbors, $i=m+1,\cdots , N$, $k \in N_+$.

In order to introduce the event-triggered mechanism, the following errors are defined as,
$$e_{{v_i}}^2 = \bar \zeta _i^2 - \zeta _i^2,{\rm{ }}t \in [\tilde t_k^i,\tilde t_{k + 1}^i),$$
where $i = m + 1, \cdots ,N$.\\

The trigger function is defined as follows,
\begin{equation} \label{type2_19}
\begin{aligned}
{\dot \phi _{i3}} = & - {\beta _{i3}}{\phi _{i3}} - {\gamma _{i3}}(e_{{v_i}}^2{^T}{H_{i2}}e_{{v_i}}^2)
+ \frac{1}{{{\theta _{i3}}}}||{\psi _{i3}}|{|^2}
 + {\tau _{i3}}{e^{ - \frac{{{\delta _{i3}}}}{{{\sigma _{i3}} + t}}}} + \frac{{{{\hat \tau }_{i3}}}}{{{{\hat \sigma}_{i3}} + t}},
\end{aligned}
\end{equation}
where ${\psi _{i3}} = \sum\limits_{j \in {\Lambda _2}} {(\bar \zeta _j^2 - \bar \zeta _i^2)}$, ${\beta _{i3}} > 0$, ${\gamma _{i3}} > 0$, ${\theta _{i3}} > 0$, ${\tau _{i3}} > 0$, ${\delta_{i3}} > 0$, ${\sigma_{i3}} > 0$, ${\hat\tau_{i3}}>0$ and ${{\hat \sigma}_{i3}}>0$ are the parameters to be specified; ${H_{i2}}$ is a constant matrix with appropriate dimension, and ${\phi _{i3}}(0) > 0$.
Then the first triggering time is set as $\tilde t_1^i = 0$. With appropriate trigger mechanism, agent $i$ will determine the triggering time sequence $\{ \tilde l_k^i\} _{k = 2}^\infty $ by,
\begin{equation} \label{type2_19.1}
\begin{aligned}
\tilde l_{k + 1}^i =& \mathop {\inf }\limits_{r \ge t_k^i} \{ \tilde r:{\omega _{i3}}({\gamma _{i3}}(e_{{v_i}}^{2T}{H_{i2}}e_{{v_i}}^2)
 - \frac{1}{{{\theta _{i3}}}}||{\psi _{i3}}|{|^2} - {\tau _{i3}}{e^{ - \frac{{{\delta _{i3}}}}{{{\sigma _{i3}} + t}}}} - \frac{{{{\hat \tau }_{i3}}}}{{{{\hat \sigma}_{i3}} + t}}) \ge {\phi _{i3}}\},
\end{aligned}
\end{equation}
where ${\omega _{i3}} > 0$ is the convergence function to be designed. \\

Let ${\beta _{i3}}{\rm{(}}t{\rm{) + }}\frac{1}{{{\omega _{i3}}}} > 0$, then ${\phi _{i3}} \ge {\phi _{i3}}(0){e^{ - \left( {{\beta _{i3}}{\rm{ + }}\frac{1}{{{\omega _{i3}}}}} \right)}}$.

The controllers of the $2^{nd}$ type followers can be designed as,
\begin{equation} \label{type2_20}
\begin{aligned}
u_i^2 =& - \hat u_i^{{a_2}} + {\hat T_{1i}}\hat x_i^2 + {\hat T_{2i}}\zeta _i^2
+ \hat M_i^3(\hat x_i^2 - {\hat X_i}\zeta _i^2) + \hat N_i^3(\sum\limits_{j \in {\Lambda _1}} {L_{ij}^2(\bar \zeta _j^2 - \bar \zeta _i^2)}),
\end{aligned}
\end{equation}
where ${\hat T_{1i}}$, ${\hat T_{2i}}$, $\hat M_i^3$ and $\hat N_i^3$ are the constant matrices with appropriate dimensions.

Similar to Lemma 3, the following lemma can be obtained.

\textbf{Lemma\,8:} According to Assumptions 4 and  5, an appropriate ${\hat T_{2i}}$ is set so that the real part of the eigenvalue of the matrix $(A_i^2 + B_i^2{\hat T_{1i}})$ is less than zero. Let ${\hat T_{2i}} = {\hat Y_i} - {\hat T_{1i}}{\hat X_i}$, $({{\hat X}_i},{{\hat Y}_i})$ is the solution of the following equations,
\begin{equation}
\begin{array}{l}
A_i^2{{\hat X}_i} + B_i^2{{\hat Y}_i} = {{\hat X}_i}{A_0},
{C_i}{{\hat X}_i} = {C_0},{\rm{   }}i = m + 1, \cdots ,N.
\end{array}
\end{equation}

\textbf{Lemma\,9:} According to Eq.(\ref{type2_20}) and Assumption 4, $\hat N_i^3$ is a matrix with appropriate dimension, which can be set to satisfy ${\hat X_i}M_i^{21}M_i^{22} = B_i^2\hat N_i^3$.\\

Let
\begin{equation}
\begin{array}{l}
{{\tilde x}^2} = {[\tilde x_{m + 1}^{2T}, \cdots ,\tilde x_N^{2T}]^T},
{\varepsilon ^2} = {[\varepsilon _{m + 1}^{2T}, \cdots ,\varepsilon _N^{2T}]^T},
\end{array}
\end{equation}

\noindent where $\tilde x_i^2 = \zeta _i^2 - {x_0}$ and $\varepsilon _i^2 = \hat x_i^2 - {\hat X_i}\zeta _i^2$, $i = m + 1, \cdots ,N$.

From Eqs.(\ref{leader_nof}), (\ref{leader_nof3}) and (\ref{type2_19}), it has,
\begin{equation} \label{type2_21}
\begin{aligned}
{{\dot {\tilde x}}^2} =& (\tilde A - ({I_{N - m}} \otimes {A_0})){\zeta ^2} + ({I_{N - m}} \otimes {A_0}){{\tilde x}^2}\\
& - {M^{21}}{M^{22}}(L_D^2 \otimes {I_n}){{\tilde x}^2} - {M^{21}}{M^{23}}({D^{21}} \otimes {I_n}){{\tilde x}^2}\\
& - {M^{21}}{M^{23}}({D^{21}} \otimes {I_n})e_v^2 + {M^{21}}{M^{23}}({A^{21}} \otimes {I_n})(e_v^1 + {{\tilde x}^1})\\
& - \left( {{1_{(N - m)}} \otimes {K_1}\sum\limits_{j \in {N_L}} {{{\bar g}_{j0}}\left( {{x_j} - {X_j}{x_0}} \right)} } \right),\\
{{\dot \varepsilon }^2} =& ({A^2} + {B^2}{{\hat T}_1}){{\hat x}^2} + ({B^2}{{\hat T}_2} - \hat X({I_{(N - m)}} \otimes {A_0})){\zeta ^2} + {W^2}({x^2} - {{\hat x}^2})\\
& + \hat X(({I_{(N - m)}} \otimes {A_0}) - \tilde A){\zeta ^2} + {B^2}{{\hat M}^3}({{\hat x}^2} - \hat X{\zeta ^2})\\
& - \hat X{M^{21}}{M^{23}}\left[ \begin{array}{l}
\sum\limits_{j = 1}^m {L_{(m + 1)j}^{21}(\zeta _j^1 - \zeta _{(m + 1)}^2)} \\
{\rm{                  }} \vdots \\
\sum\limits_{j = 1}^m {L_{Nj}^{21}(\zeta _j^1 - \zeta _N^2)}
\end{array} \right],
\end{aligned}
\end{equation}
where $\tilde A = diag\{ {\tilde A_{m + 1}}, \cdots ,{\tilde A_N}\}$, ${A^2} = diag\{ A_{m + 1}^2, \cdots ,A_N^2\}$, ${B^2} = diag\{ B_{m + 1}^2, \cdots ,B_N^2\}$, \\
${M^{21}} = diag\{ M_{m + 1}^{21}, \cdots ,M_N^{21}\}$, ${M^{22}} = diag{\rm{\{ }}M_{m + 1}^{22}, \cdots ,M_N^{22}\}$, ${M^{23}} = diag{\rm{\{ }}M_{m + 1}^{23}, \cdots ,M_N^{23}\}$,\\ 
${\hat T_1} = diag\{ {\hat T_{1(m + 1)}}, \cdots ,{\hat T_{1N}}\}$ and ${\hat T_2} = diag\{ {\hat T_{2(m + 1)}}, \cdots ,{\hat T_{2N}}\}$.

% \noindent \textbf{Remark\,XX:}

Considering the following Lyapunov functions,
\begin{equation} \label{type2_21.1}
\begin{array}{l}
{\hat V_{2i}} = \hat V_{2i}^1 + \hat V_{2i}^2,
{\rm{  }}\hat V_{2i}^1 = \tilde x_i^{2T}{\hat P_{1i}}{\tilde x^2} + {\phi _{i3}},
{\rm{  }}\hat V_{2i}^2 = \varepsilon _i^{2T}{\hat P_{2i}}\varepsilon _i^2.
\end{array}
\end{equation}

Let ${\hat V_2} = \sum\limits_{j = 1}^m {{{\hat V}_{2j}}}$, $\hat V_2^1 = \sum\limits_{j = 1}^m {\hat V_{2j}^1}$, $\hat V_2^2 = \sum\limits_{j = 1}^m {\hat V_{2j}^2}$.

The time derivative of ${\hat V}_2^1$ along with System (\ref{type2_21}) is given by
\begin{footnotesize}
\begin{equation}
\begin{array}{l}
\dot {\hat V}_2^1 = 2{\zeta ^{2T}}{(\tilde A - ({I_{N - m}} \otimes {A_0}))^T}{{\hat P}_1}{{\tilde x}^2}
+ {{\tilde x}^{2T}}(F_1^T{{\hat P}_1} + {{\hat P}_1}{F_1}){{\tilde x}^2} \\
- 2{{\tilde x}^{2T}}{{\hat P}_1}{{\hat Q}_1}e_v^2 
+ 2{{\tilde x}^{2T}}{M^{21}}{M^{23}}({A^{21}} \otimes {I_n})(e_v^1 + {{\tilde x}^1}) - 2{{\tilde x}^{2T}}{K^1}\Xi \\
\sum\limits_{i = m + 1}^N {\left( { - {\beta _{i3}}{\phi _{i3}} - {\gamma _{i3}}(e_{{v_i}}^2{{}^T}{H_{i2}}e_{{v_i}}^2) + \frac{1}{{{\theta _{i3}}}}||{\psi _{i3}}|{|^2} + {\hat \delta_{i3}} } \right)},
\end{array}
\end{equation}
\end{footnotesize}
where ${\hat Q_1} = {M^{21}}{M^{23}}({D^{21}} \otimes {I_n}), {F_1} = ({I_{N - m}} \otimes {A_0}) - {M^{21}}{M^{22}}(L_D^2 \otimes {I_n}) - {M^{21}}{M^{23}}({D^{21}} \otimes {I_n})$, ${\hat \delta_{i3}} = {\tau _{i3}}{e^{ - \frac{{{\delta _{i3}}}}{{{\sigma _{i3}} + t}}}} + \frac{{{{\hat \tau }_{i3}}}}{{{{\hat \sigma}_{i3}} + t}}$.

Let\\
${\hat \theta _{\min }} = \min \{ {\theta _{(m + 1)3}}, \cdots ,{\theta _{N3}}\}$,
${H_2} = diag({H_{(m + 1)2}}, \cdots ,{H_{N2}})$,
$ {\hat \gamma _{\min }} = \min \{ {\gamma _{(N - m)3}}, \cdots ,{\gamma _{(N - m)3}}\}$.

Then,
\begin{equation}
\begin{array}{l}
\sum\limits_{i = m + 1}^N {\frac{1}{{{\theta _{i3}}}}||{\psi _{i3}}|{|^2}}  \le \frac{1}{{{{\hat \theta }_{\min }}}}\sum\limits_{i = m + 1}^N {||{\psi _{i3}}|{|^2}}\\
 = \frac{1}{{{{\hat \theta }_{\min }}}}\hat \chi _1^T{[{I_{N - m}},{I_{N - m}}]^T}({L^2}{L^2} \otimes {I_{N - m}})[{I_{N - m}},{I_{N - m}}]{{\hat \chi }_1},
\end{array}
\end{equation}
where ${\hat \chi _1} = {[{\tilde x^{2T}},e_v^{2T}]^T}$. Then it has,
\begin{equation}
\begin{array}{l}
{{\tilde x}^{2T}}(F_1^T{{\hat P}_1} + {{\hat P}_1}{F_1}){{\tilde x}^2} - 2{{\tilde x}^{2T}}{{\hat P}_1}{{\hat Q}_1}e_v^2
 - \sum\limits_{j = m + 1}^N {({\gamma _{i3}}e_{{v_i}}^2{{}^T}{H_{i2}}e_{{v_i}}^3)}  \le \hat \chi _1^T{{\hat \Gamma }^1}\hat \chi _1^T,
\end{array}
\end{equation}
where ${\hat \Gamma ^1} = [\hat \Gamma _{11}^1,*;\hat \Gamma _{21}^1,\hat \Gamma _{22}^1]$ with $\hat \Gamma _{11}^1 = F_1^T{\hat P_1} + {\hat P_1}{F_1} + {c_2}{\hat P_1}$, $\hat \Gamma _{21}^1 =  - {\hat Q_1}{\hat P_1}$, $\hat \Gamma _{22}^1 =  - {\hat \gamma _{\min }}{H_2}$, ${c_2} > 0$.

By Young's inequality,
\begin{equation}\label{type2_22}
\begin{array}{l}
2{{\tilde x}^{2T}}{M^{21}}{M^{23}}({A^{21}} \otimes {I_n})e_v^1
 \le {{\hat \beta }_1}{{\tilde x}^{2T}}{{\tilde x}^T} + \frac{1}{{{{\hat \beta }_1}}}||{M^{21}}{M^{23}}({A^{21}} \otimes {I_n})e_v^1|{|^2},
\end{array}
\end{equation}
\begin{equation}\label{type2_23}
\begin{array}{l}
2{{\tilde x}^{2T}}{M^{21}}{M^{23}}({A^{21}} \otimes {I_n}){{\tilde x}^1}
 \le {{\hat \beta }_2}{{\tilde x}^{2T}}{{\tilde x}^2} + \frac{1}{{{{\hat \beta }_2}}}||{M^{21}}{M^{23}}({A^{21}} \otimes {I_n}){{\tilde x}^1}|{|^2},
\end{array}
\end{equation}
where ${\hat \beta _1}$ and ${\hat \beta _2}$ are the real numbers to be determined.

According to Eqs.(\ref{type1_18}), (\ref{type1_19}), (\ref{type2_22}), (\ref{type2_23}) and Lemma 5, it has,
\begin{equation}
\begin{aligned}
\dot {\hat V}_2^1\le&  - {\mu _2}{c_2}{{\tilde x}^{2T}}{{\hat P}_1}{{\tilde x}^2} + {\beta _3}||{{\tilde x}^{1T}}{K^1}|{|^2}\\
 &+ {\beta _4}||{{\tilde x}^{1T}}{K^1}|{|^2} + \frac{1}{{{{\hat \beta }_2}}}||{M^{21}}{M^{23}}({A^{21}} \otimes {I_n}){{\tilde x}^1}|{|^2} + {\hat \Omega _1}\\
 &+ {\hat \Omega _2} + \frac{1}{{{\beta _4}}}||{{\bar g}_1}{{\hat K}_1}{\varepsilon ^1}|{|^2} -
  \sum\limits_{i = m + 1}^N {\left( {{\beta _{i3}}{\phi _{i3}} - {\hat \delta_{i3}}} \right)},
\end{aligned}
\end{equation}
where ${\mu _2} = 1 - {\mu _1}$, ${\mu _1} > 0, {\hat \Omega _1} = \frac{1}{{{{\hat \beta }_1}}}||{M^{21}}{M^{23}}({A^{21}} \otimes {I_n})e_v^1|{|^2} + \frac{1}{{{\beta _3}}}||{\bar g_1}{\hat K_1}{\bar x^1}|{|^2} + 2{\zeta ^{2T}}{(\tilde A - ({I_{N - m}} \otimes {A_0}))^T}{\hat P_1}{\tilde x^2},\\ {\hat \Omega _2} =  - {\mu _1}{c_2}{\tilde x^{2T}}{\hat P_1}{\tilde x^2} + {\hat \beta _1}{\tilde x^{2T}}{\tilde x^2} + {\hat \beta _2}{\tilde x^{2T}}{\tilde x^2}$.

According to Lemma 1, Lemma 7 and Eq.(\ref{type1_10.2}) ${\hat \Omega _1}$ is bounded.

Then it has,
\begin{equation}\label{type2_24}
\begin{aligned}
\dot {\hat V}_2^1 \le &  - {\mu _1}{c_2}{{\tilde x}^{2T}}{{\hat P}_1}{{\tilde x}^2} + {\beta _3}||{{\tilde x}^{1T}}{K^1}|{|^2} + {\beta _4}||{{\tilde x}^{1T}}{K^1}|| \\
&+ \frac{1}{{{\beta _4}}}||{{\bar g}_1}{{\hat K}_1}{\varepsilon ^1}|{|^2} +\frac{1}{{{{\hat \beta }_2}}}||{M^{21}}{M^{23}}({A^{21}} \otimes {I_n}){{\tilde x}^1}|{|^2}
+ {{\hat \Omega }_2}\\
& + {{\hat \Omega }_1} - \sum\limits_{i = m + 1}^N {\left( {{\beta _{i3}}{\phi _{i3}} - {\hat \delta_{i3}}} \right)}.
\end{aligned}
\end{equation}

\textbf{Lemma\,10:} According to Assumption 1, Assumption 4, Definition 1 and Lemma 8, an appropriate matrix ${\hat M^3}$ can be found so that there is a positive definite matrix ${\hat P_2}$, it has,
\begin{small}
\begin{equation}
\begin{array}{l}
{\left( {{A^2} + {B^2}({{\hat T}_1} + {{\hat M}^3})} \right)^T}{\hat P_2} + {\hat P_2}\left( {{A^2} + {B^2}({{\hat T}_1} + {{\hat M}^3})} \right) 
\le - {\tilde \mu _1}{\hat P_2},
\end{array}
\end{equation}
\end{small}
where ${\tilde \mu _1} > 0$. According to Lemma 8, the time derivative of ${\hat V}_2^2$ along the trajectory of System (\ref{type2_21}) is calculated as,
\begin{equation}
\begin{aligned}
\dot {\hat V}_2^2 =& {\varepsilon ^2}^T{\left( {{A^2} + {B^2}({T_1} + {{\hat M}^3})} \right)^T}{{\hat P}_2}{\varepsilon ^2} + 2{\varepsilon ^2}^T{{\hat P}_2}{W^2}({x^2} - {{\hat x}^2})\\
& + {\varepsilon ^2}^T{{\hat P}_2}\left( {{A^2} + {B^2}({{\hat T}_1} + {{\hat M}^3})} \right){\varepsilon ^2}
+ 2{\varepsilon ^2}^T{{\hat P}_2}\hat X(({I_{(N - m)}} \otimes {A_0}) - \tilde A){\zeta ^2}\\
& - 2{\varepsilon ^2}^T\hat X{M^{21}}{M^{23}}\left[ \begin{array}{l}
\sum\limits_{j = 1}^m {L_{(m + 1)j}^{21}(\bar \zeta _j^1 - \bar \zeta _{(m + 1)}^2)} \\
{\rm{                 }} \vdots \\
\sum\limits_{j = 1}^m {L_{Nj}^{21}(\zeta _j^1 - \zeta _N^2)}
\end{array} \right].
\end{aligned}
\end{equation}

By Young's inequality and Lemma 10, it has,
\begin{equation} \label{type2_25}
\begin{aligned}
\dot {\hat V}_2^2 \le & - {{\tilde \mu }_1}{\varepsilon ^2}^T{{\hat P}_2}{\varepsilon ^2} + {{\hat \tau }_6}{\varepsilon ^2}^T{\varepsilon ^2}
 + \tilde M + \frac{1}{{{{\hat \tau }_2}}}||{{\hat P}_2}\hat X{M^{21}}{M^{23}}({A^{21}} \otimes {I_{(N - m)}}){{\tilde x}^1}|{|^2}\\
& + \frac{1}{{{{\hat \tau }_3}}}||{{\hat P}_2}\hat X{M^{21}}{M^{23}}({A^{21}} \otimes {I_{(N - m)}}){{\tilde x}^2}|{|^2},
\end{aligned}
\end{equation}
where\\
$\tilde M = \frac{1}{{{{\hat \tau }_5}}}||{\hat P_2}{W^2}{\bar x^2}|{|^2} + \frac{1}{{{{\hat \tau }_1}}}||{\hat P_2}\hat X{M^{21}}{M^{23}}({A^{21}} \otimes {I_{(N - m)}})({\bar \zeta ^1} - {\zeta ^1})|{|^2} + \frac{1}{{{{\hat \tau }_4}}}||{\hat P_2}\hat X{M^{21}}{M^{23}}({A^{21}} \otimes {I_{(N - m)}})({\zeta ^2} - {\bar \zeta ^2})|{|^2} + 2{\varepsilon ^2}^T{\hat P_2}\hat X(({I_{(N - m)}} \otimes {A_0}) - \tilde A){\zeta ^2}, {\hat \tau _6} = {\hat \tau _1} + {\hat \tau _2} + {\hat \tau _3} + {\hat \tau _4} + {\hat \tau _5}$.

According to Lemma 2, Lemma 7, Eq.(\ref{type1_10.2}) and Eq.(\ref{type2_19.1}), $\tilde M$ is bounded.

From Eqs.(\ref{type2_21.1}), (\ref{type2_24}) and (\ref{type2_25}), it has,
\begin{equation} \label{type2_26}
\begin{aligned}
{{\hat V}_2} =& \hat V_2^1 + \hat V_2^2
 \le  - {\mu _2}{c_2}{{\tilde x}^{2T}}{{\hat P}_1}{{\tilde x}^2} + {\beta _3}||{{\tilde x}^{1T}}{K^1}|{|^2} + {{\hat \Omega }_1}\\
& + \frac{1}{{{{\hat \beta }_2}}}||{M^{21}}{M^{23}}({A^{21}} \otimes {I_n}){{\tilde x}^1}|{|^2}  + \frac{1}{{{\beta _4}}}||{{\bar g}_1}{{\hat K}_1}{\varepsilon ^1}|{|^2} + {{\hat \Omega }_2}\\
& + \frac{1}{{{{\hat \tau }_2}}}||{{\hat P}_2}\hat X{M^{21}}{M^{23}}({A^{21}} \otimes {I_{(N - m)}}){{\tilde x}^1}|{|^2} - {{\tilde \mu }_1}{\varepsilon ^2}^T{{\hat P}_2}{\varepsilon ^2}\\
& + \frac{1}{{{{\hat \tau }_3}}}||{{\hat P}_2}\hat X{M^{21}}{M^{23}}({A^{21}} \otimes {I_{(N - m)}}){{\tilde x}^2}|{|^2} + {\beta _4}||{{\tilde x}^{1T}}{K^1}||\\
& + {{\hat \tau }_6}{\varepsilon ^2}^T{\varepsilon ^2} + \tilde M - \sum\limits_{i = m + 1}^N {\left( {{\beta _{i3}}{\phi _{i3}} - {\hat \delta_{i3}}} \right)}.
\end{aligned}
\end{equation}

\textbf{Lemma\,11:} The controller designed in Eq.(\ref{type1_11}) and Eq.(\ref{type2_20}) can ensure that $(\tilde x_i^1,\varepsilon _i^1,\tilde x_j^2,\varepsilon _j^2)$, $i = 1, \cdots ,m$, $j = m + 1, \cdots ,N$ is UUB under Assumptions 1-4, Eq.(\ref{type1_10.2}) and Eq.(\ref{type2_19.1}).

\textbf{Proof\,:} See Appendix C.

Next, it is proved that the Zeno phenomenon will not occur in the heterogeneous multi-agent system with the proposed event triggered mechanism on the premise of system consistency. So the following lemma is applied.

\textbf{Lemma\,12:} There is no Zeno behavior for each agent under the event-triggered mechanism (\ref{type1_10.2}).

\textbf{Proof\,:} See Appendix D.

\textbf{Theorem\,1:} Under Assumptions 1-5, the control protocol in Eq.(\ref{type1_11}) and Eq.(\ref{type2_20}) can make the followers System (\ref{leader_nof}) and leader System (\ref{leader_nof3}) under the communication failure (\ref{leader_1}) and actuator failure (\ref{leader_2}) to reach the consistency with UUB.

\textbf{Proof\,:} First, the $1^{st}$ type followers are discussed. From System (\ref{leader_nof}) and System (\ref{leader_nof3}), it has,
\begin{small}
\begin{equation}
\begin{aligned}
||y_i^1 - {y_0}|{|^2}
 =& ||{C_i}x_i^1 - {C_i}\hat x_i^1 + {C_i}\hat x_i^1 - {C_i}{X_i}\zeta _i^1
 + {C_i}{X_i}\zeta _i^1 - {C_0}{x_0}|{|^2}\\
 \le& 3||{C_i}|{|^2}||x_i^1 - \hat x_i^1|{|^2} + 3||{C_i}|{|^2}||\hat x_i^1 - {X_i}\zeta _i^1|{|^2}
 + 3||{C_i}{X_i}|{|^2}||\zeta _i^1 - {x_0}|{|^2}.
\end{aligned}
\end{equation}
\end{small}

According to Lemma 1 and Lemma 11, let ${T_2} = \max \{ {T_i},T\}$, when $t > {T_2}$, $||y_i^1 - {y_0}|{|^2}$ is UUB, $i = 1, \cdots ,m$. Similarly, it can be shown that $||y_j^2 - {y_0}|{|^2}$, is UUB, $j = m + 1, \cdots ,N$. This demonstrates that all the followers and the leader can achieve consistency with UUB.      $\hfill{} \Box$

\section{Simulation experiments and result analysis}
To illustrate the performance of the proposed distributed protocol, simulation experiments are performed in this Section. In particular, consider a multi-agent system composed of one leader, three $1^{st}$ type followers and five $2^{nd}$ type followers, that is, ${\Lambda _1} = \{ 1,2,3\}$, ${\Lambda _2} = \{ 4,5,6,7,8\}$, the communication topology graph for such multi-agent system can refer in Fig.\ref{fig_1}. Then a heterogeneous multi-agent system with the dynamics matrices, ${A_i} = [{a_{i1}},{a_{i2}};{a_{i3}},{a_{i4}}]$, ${B_i}=[{b_{i1}},{b_{i2}};{b_{i3}},{b_{i4}}]$,${C_i} = [1,0]$, $i = 1, \cdots ,8$, where $({a_{i1}},{a_{i2}},{a_{i3}},{a_{i4}},{b_{i1}},{b_{i2}},{b_{i3}},{b_{i4}})$ is chosen as,\\
$({\rm{1,  - 1,  - 2, 3}},{\rm{ 2, 0}}{\rm{.5, 0}}{\rm{.5, 1}})$, $({\rm{1}}{\rm{.6, 1}}{\rm{.3, 1}}{\rm{.1, 1}}{\rm{.5}},{\rm{ 2}}{\rm{.0, 1}}{\rm{.8, 1}}{\rm{.5, 1}}{\rm{.2}})$, $({\rm{1}}{\rm{.4 1}}{\rm{.2, 1}}{\rm{.3, 1}}{\rm{.1,1}}{\rm{.4, 1}}{\rm{.1, 1}}{\rm{.1, 1}}{\rm{.5}})$,\\ $({\rm{1}}{\rm{.0, 1}}{\rm{.5, 1}}{\rm{.3, 1}}{\rm{.6, 2}}{\rm{.4, 1}}{\rm{.7, 1}}{\rm{.7, 3}}{\rm{.0}})$, $({\rm{2}}{\rm{.5, 2, 2, 2}}{\rm{.6, 1}}{\rm{.77, 2}}{\rm{.4, 2}}{\rm{.4, 2}}{\rm{.3}})$, $({\rm{1}}{\rm{.5, 2}}{\rm{.4, 2}}{\rm{.3, 2}}{\rm{.6, 2}}{\rm{.5, 1}}{\rm{.3, 1}}{\rm{.3, 3}}{\rm{.0}})$,\\
$({\rm{1}}{\rm{.6, 2}}{\rm{.2, 1}}{\rm{.5, 2}}{\rm{.3, 2}}{\rm{.3, 2}}{\rm{.7, 2}}{\rm{.7, 2}})$, $({\rm{1}}{\rm{.7, 2}}{\rm{.6, 2}}{\rm{.2, 2}}{\rm{.8, 2}}{\rm{.3, 2}}{\rm{.6, 2}}{\rm{.6, 2}})$.

\begin{figure}[!ht]
\centering
 \includegraphics[width=0.7\textwidth]{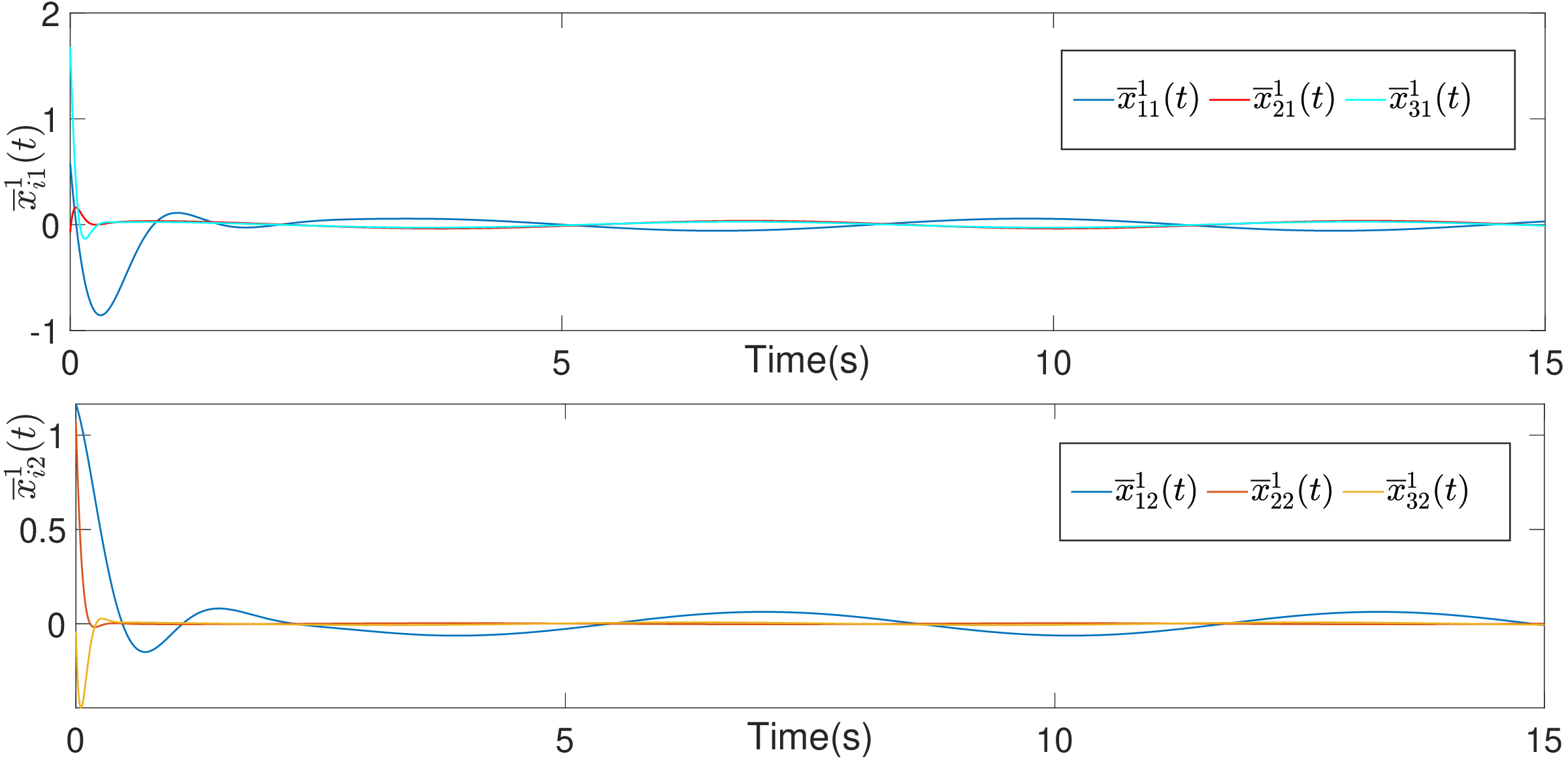}
\caption{The estimation errors $\bar x_i^1 = x_i^1 - \hat x_i^1$ of the $1^{st}$ type followers, $i = 1, \cdots ,3$.}
\label{Pi_1_x_heng}
\end{figure}

\begin{figure}[!ht]
\centering
\includegraphics[width=0.7\linewidth,height = 4cm]{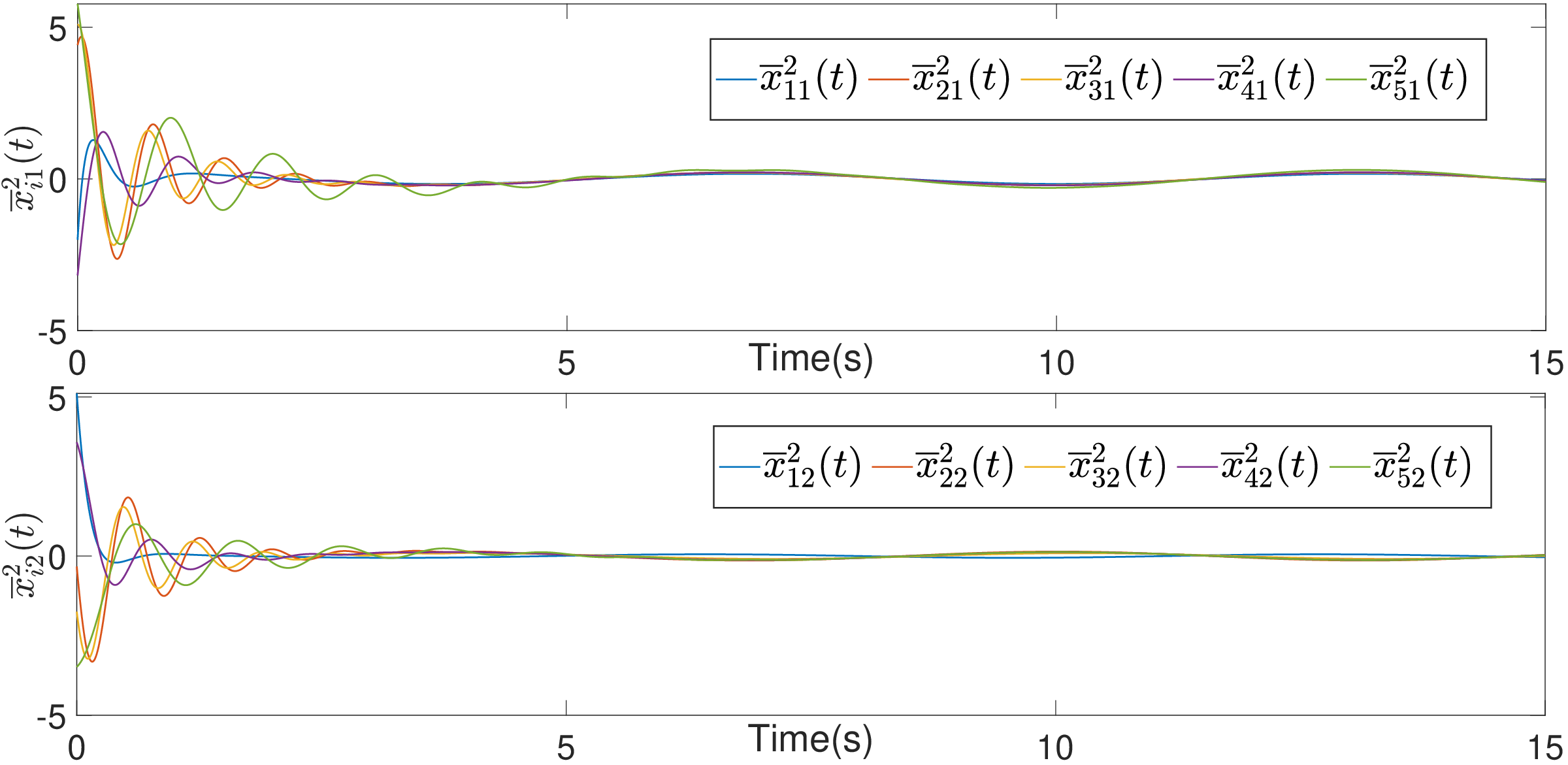}
\caption{The estimation errors $\bar x_i^2 = x_i^2 - \hat x_i^2$ of the $2^{nd}$ type followers, $i = 4, \cdots ,8$.}
\label{Pi_2_x_heng}
\end{figure}

The reference model is considered as, ${A_0} = [0,2; - 1.5,0]$, ${C_0} = [1,2]$, the communication fault is $\delta _{ij}^a = 0.25\cos r_{ij}^1t$, $\delta _{kj}^a = 0.25\sin r_{kj}^2t$, $\delta _{0i}^g = 0.25\sin r_i^3t$, $\delta _{i0}^g = 0.25\cos r_i^4t$, and the actuator fault is designed as, ${U^a} = [0.3\sin {r_5t};0.4\cos {r_6t}]$, where $r_{ij}^1$, $r_{kj}^2$, $r_i^3$, $r_i^4$, ${r_5}$, ${r_6}$ are random values between $[0,1]$, $i = 1,2,3$, $j = 1, \cdots ,8$, $k = 4, \cdots ,8$. It can easily be verified that Assumptions 3-6 are satisfied. With appropriate gain matrices selection, Lemma 1 and Lemma 2 are hold. The trajectories of the estimation errors $\bar x_i^1 = {[\bar x_{i1}^{1T},\bar x_{i2}^{1T}]^T}$ and $\bar x_j^2 = {[\bar x_{j1}^{2T},\bar x_{j2}^{2T}]^T}$ are demonstrated in Fig.\ref{Pi_1_x_heng} and Fig.\ref{Pi_2_x_heng}, where $\bar x_{ik}^1$ represents the ${k^{th}}$ dimension element of $\bar x_i^1$, and $\bar x_{jk}^{2}$ is the ${k^{th}}$ dimension element of $\bar x_{j}^{2}$, $i \in {\Lambda_1}$, $j \in {\Lambda_2}$, $k = 1, 2$.

\begin{figure}[!ht]
\centering
\includegraphics[width=0.7\linewidth,height = 4cm]{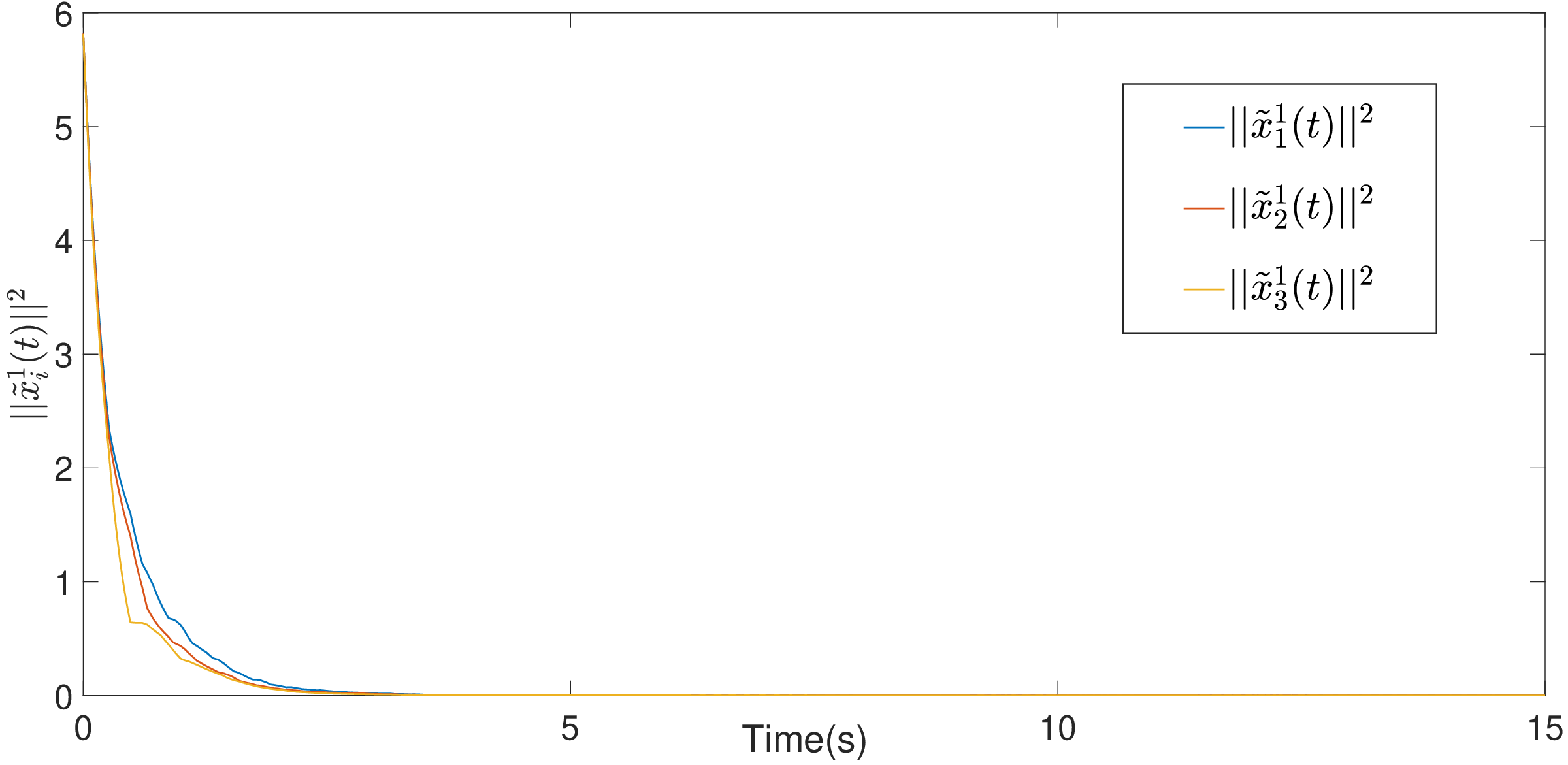}
\caption{Trajectories of observer errors $||\tilde x_i^1|{|^2}$ of the $1^{st}$ type followers, $i = 1, \cdots ,3$.}
\label{Pi_1_x_wave}
\end{figure}

\begin{figure}[!ht]
\centering
\includegraphics[width=0.7\linewidth,height = 4cm]{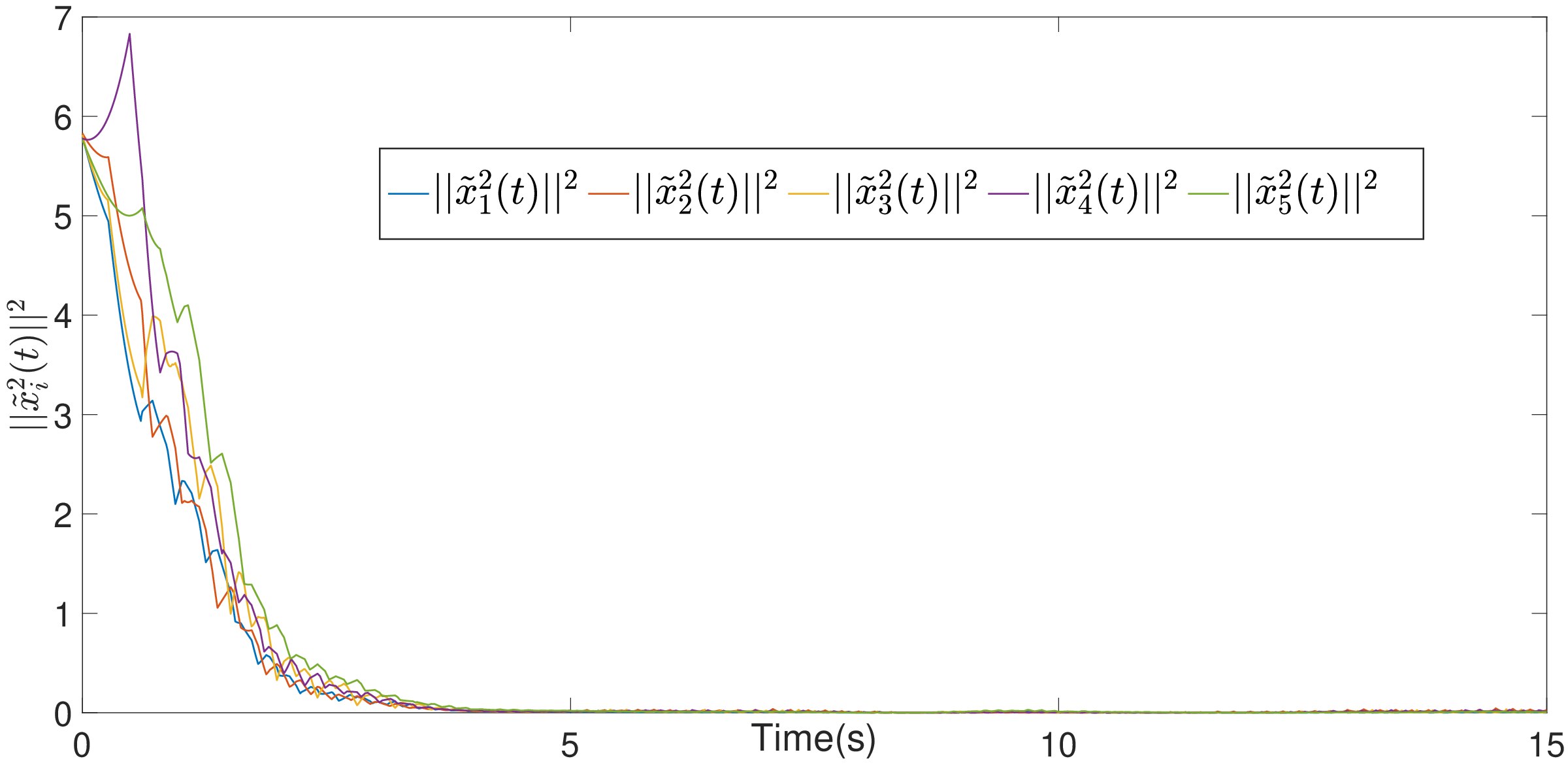}
\caption{Trajectories of the observer errors $||\tilde x_i^2|{|^2}$ of the $2^{nd}$ type followers, $i = 4, \cdots ,8$.}
\label{Pi_2_x_wave}
\end{figure}

\begin{figure}[!ht]
\centering
\includegraphics[width=0.7\linewidth,height = 4.5cm]{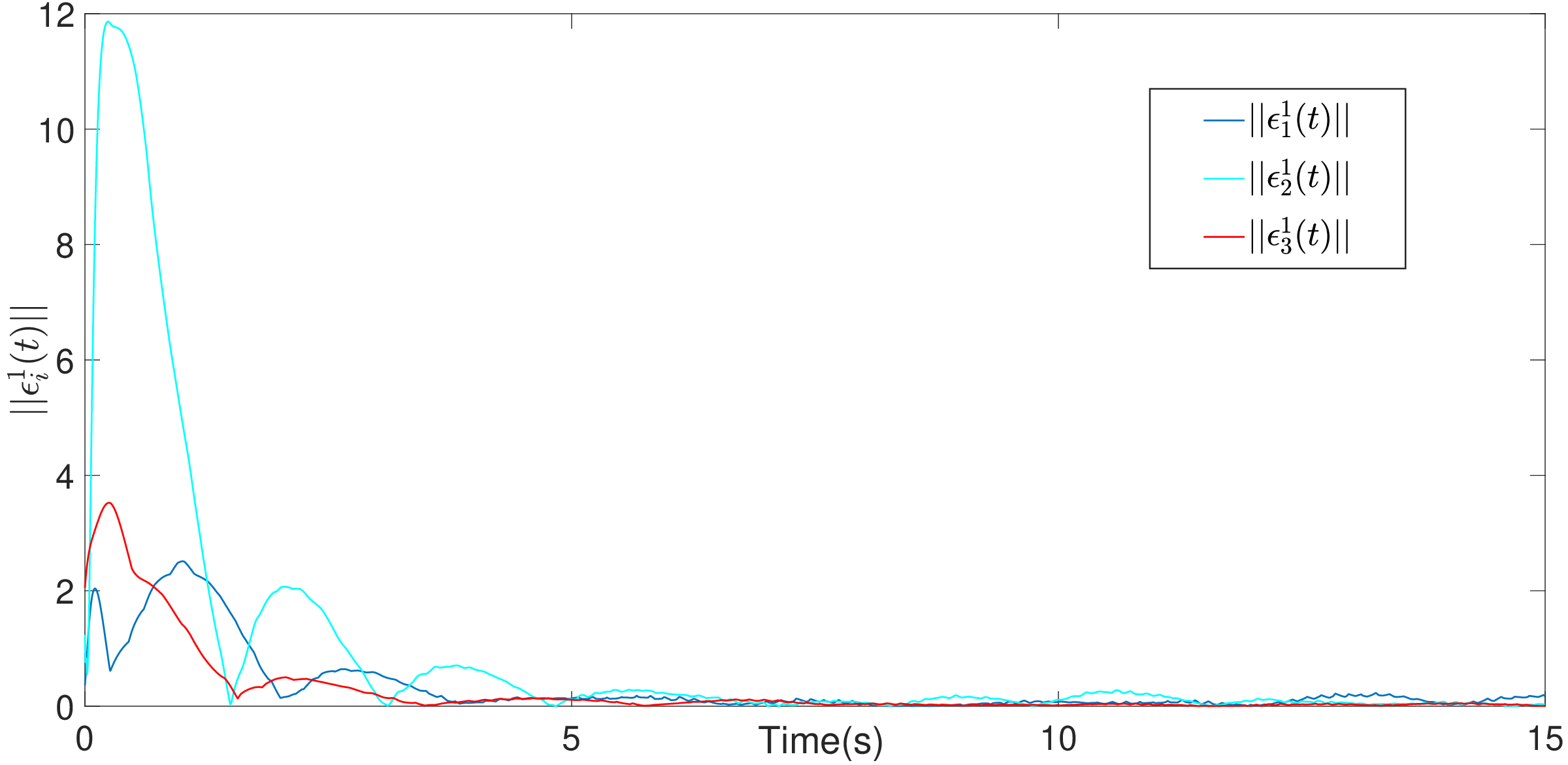}
\caption{The errors $||\varepsilon _i^1|{|^2}$ of the $1^{st}$ type followers, $i = 1, \cdots ,3$.}
\label{Pi_1_eplison}
\end{figure}

\begin{figure}[!ht]
\centering
\includegraphics[width=0.7\linewidth,height = 4.5cm]{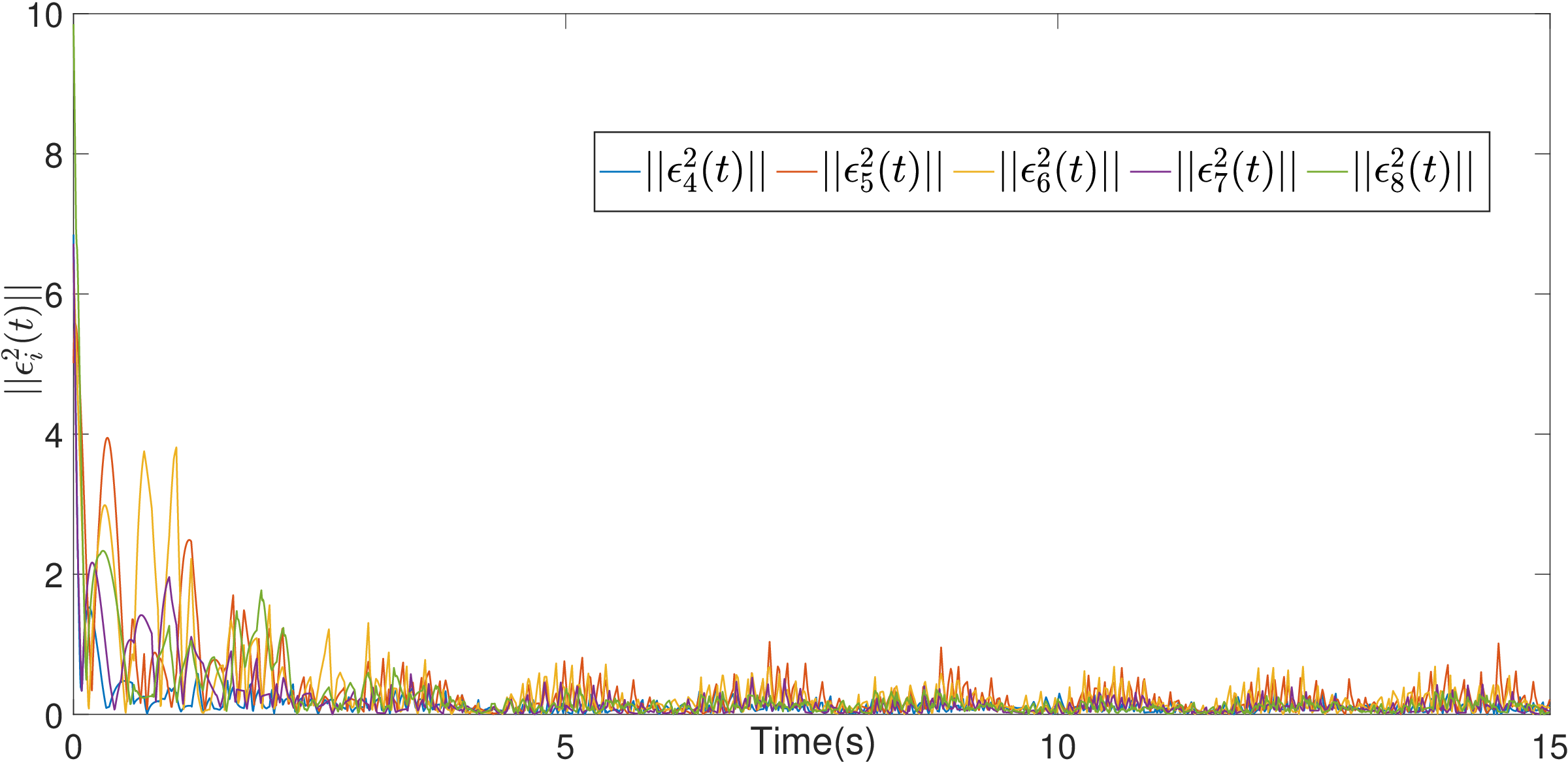}
\caption{The errors $||\varepsilon _i^2|{|^2}$ of the $2^{nd}$ type followers, $i = 4, \cdots ,8$.}
\label{Pi_2_eplison}
\end{figure}

\begin{figure}[!ht]
\centering
\includegraphics[width=0.7\linewidth,height = 4.2cm]{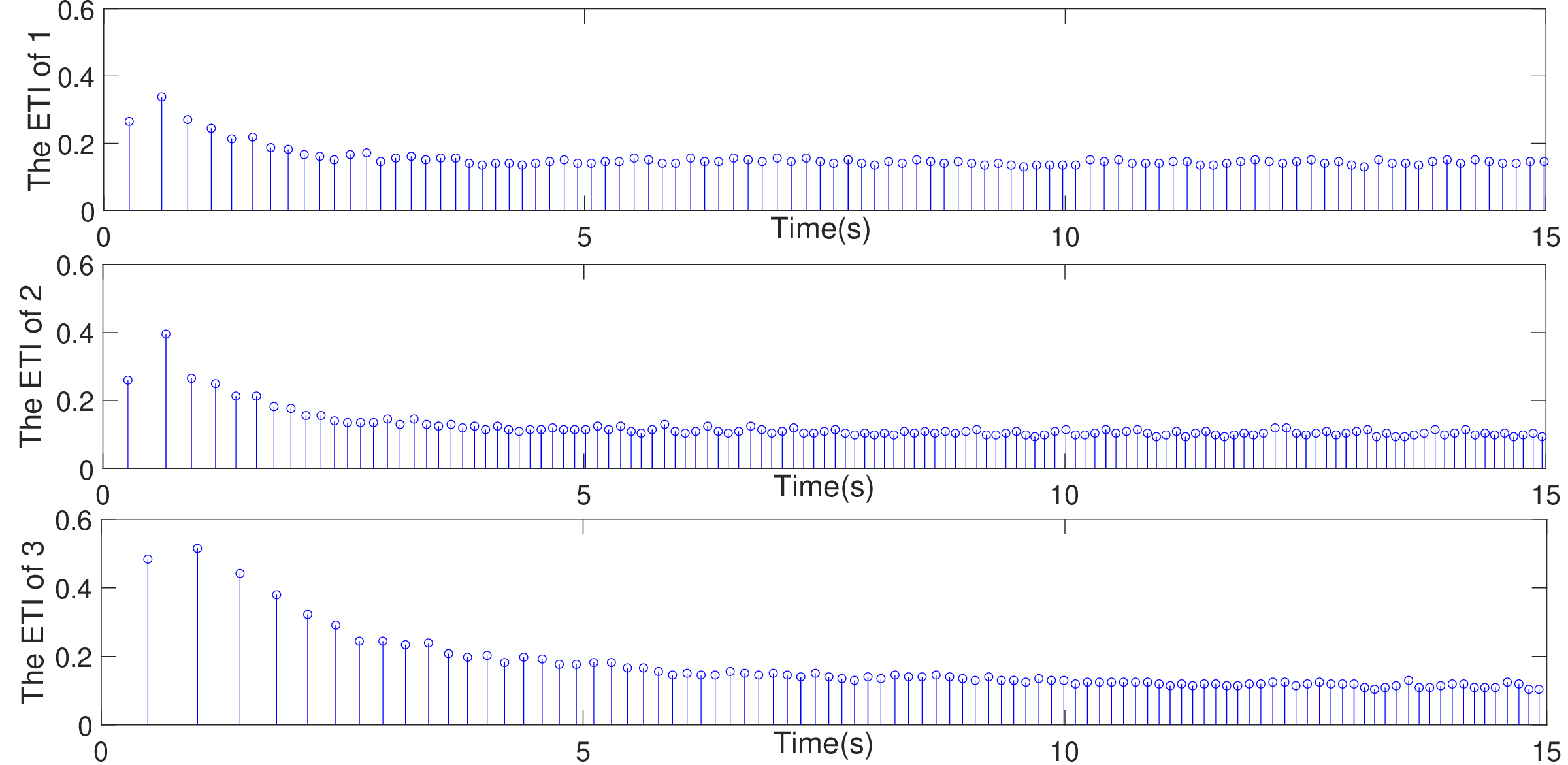}
\caption{The triggering intervals of the leader observation state $\zeta _i^1$ of the $1^{st}$ type followers, $i = 1, \cdots ,3$.}
\label{Pi1_event}
\end{figure}

\begin{figure}[!ht]
\centering
\includegraphics[width=0.7\linewidth,height = 4cm]{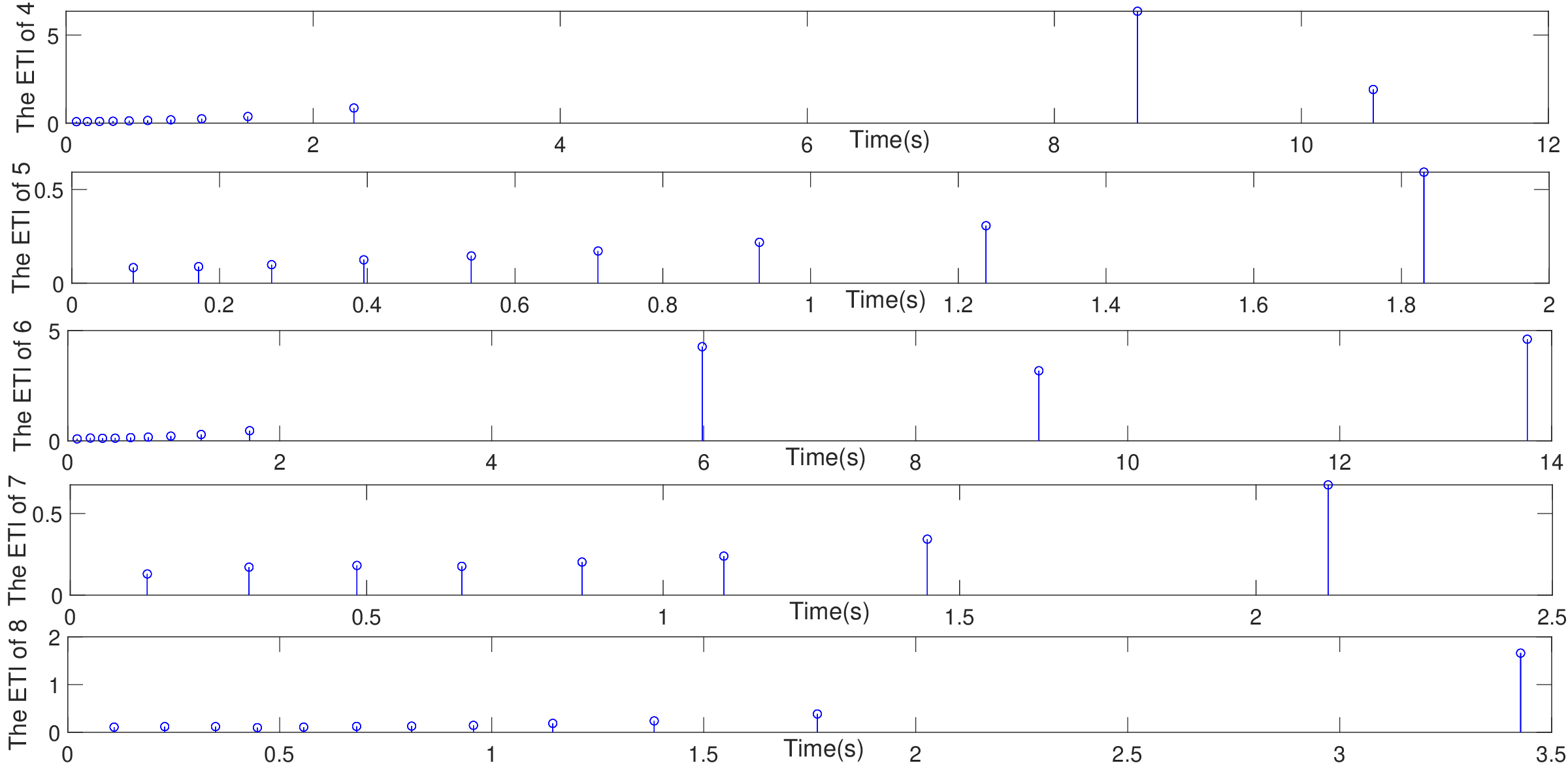}
\caption{The triggering intervals of the estimating state ${\hat A_i}$ of the $2^{nd}$ type followers, $i = 4, \cdots ,8$.}
\label{Pi2_event_A0}
\end{figure}

\begin{figure}[!ht]
\centering
\includegraphics[width=0.7\linewidth,height = 4.1cm]{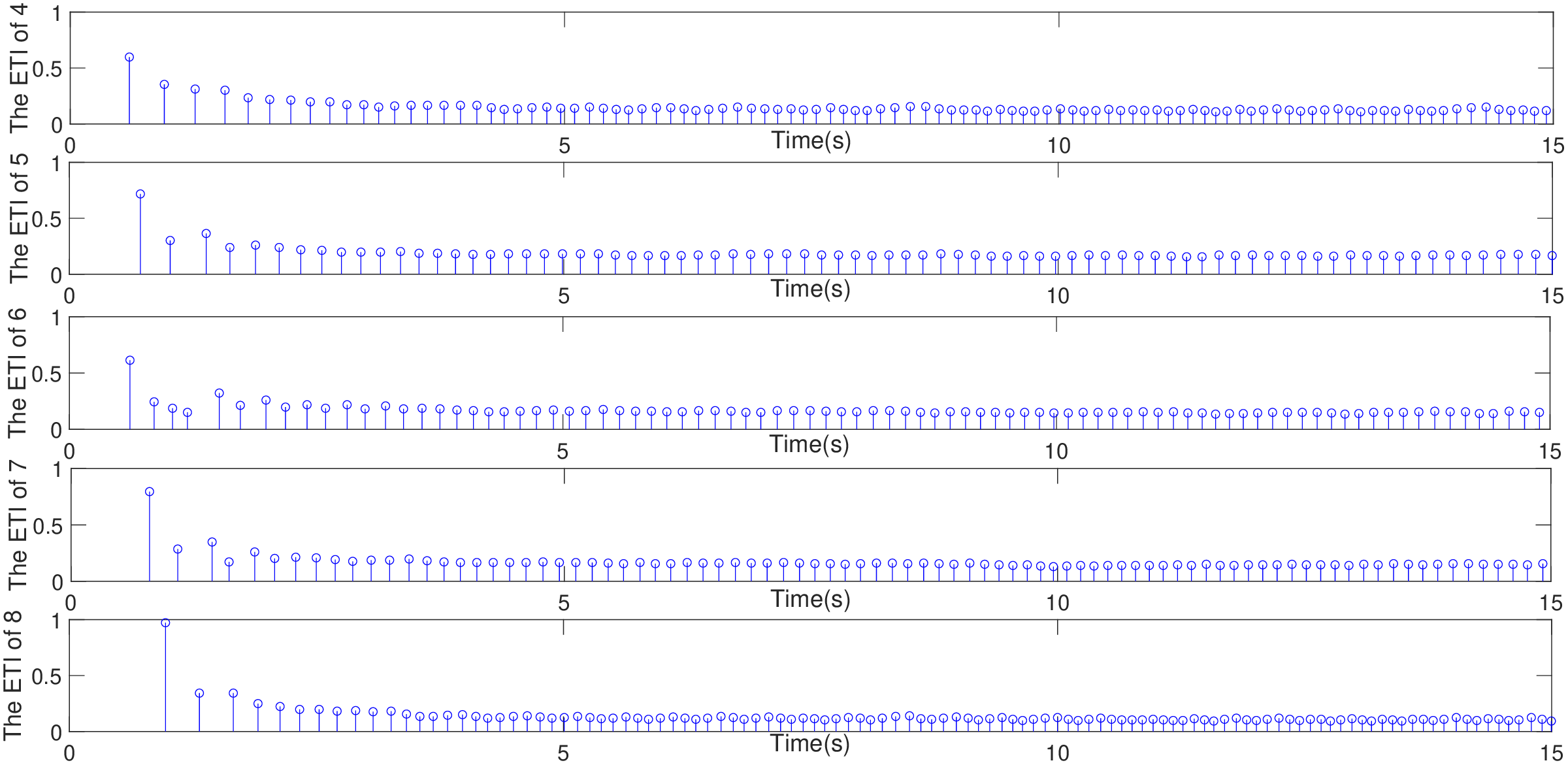}
\caption{The triggering intervals of the leader observation state $\zeta _i^2$ of the $2^{nd}$ type followers, $i = 4, \cdots ,8$.}
\label{Pi2_event}
\end{figure}

\begin{figure}[!ht]
\centering
\includegraphics[width=0.7\linewidth,height = 4.5cm]{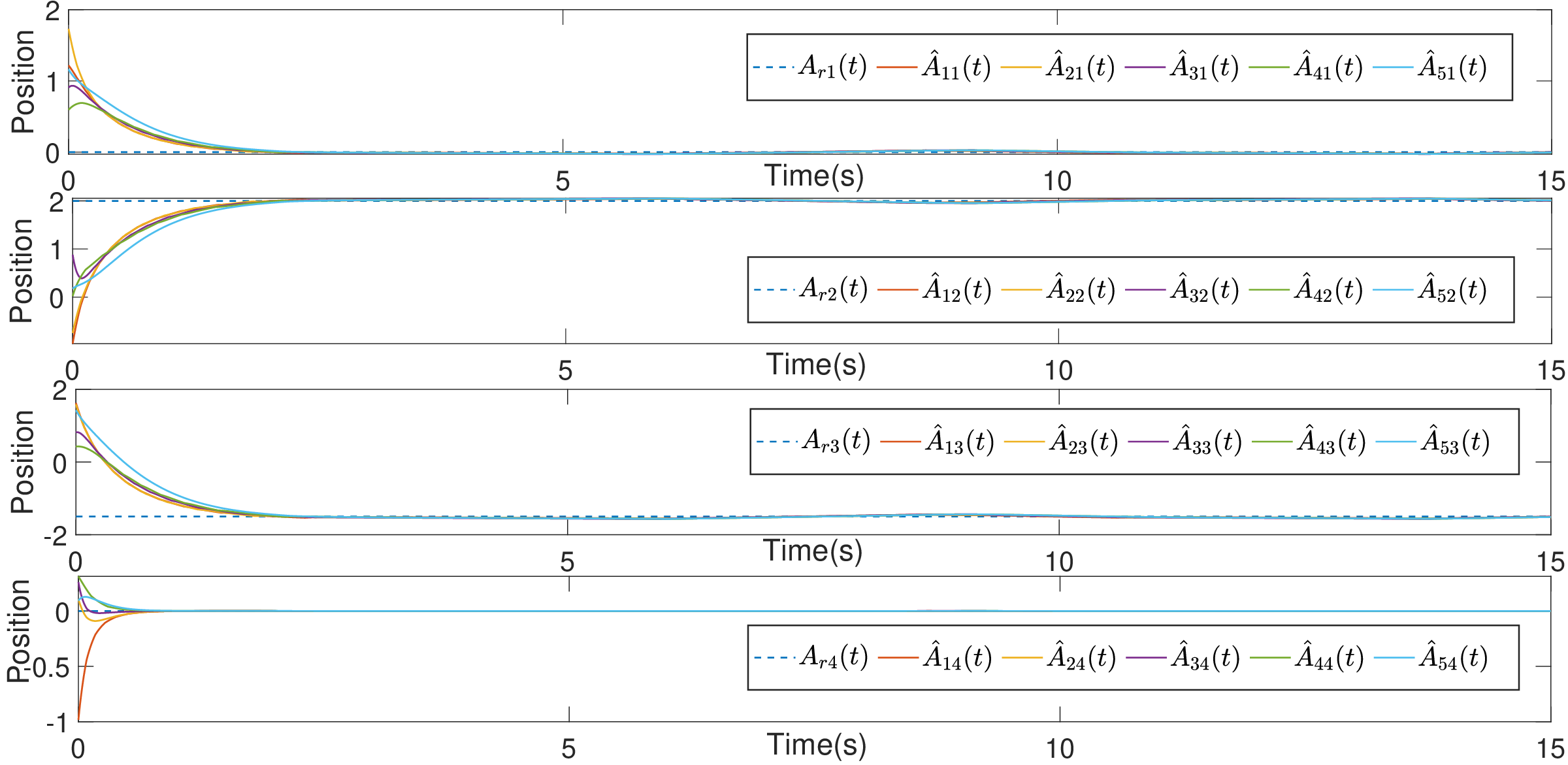}
\caption{The estimation of ${A_r}$.}
\label{Pi_2_A0}
\end{figure}

\begin{figure}[!ht]
\centering
\includegraphics[width=0.7\linewidth,height = 4.5cm]{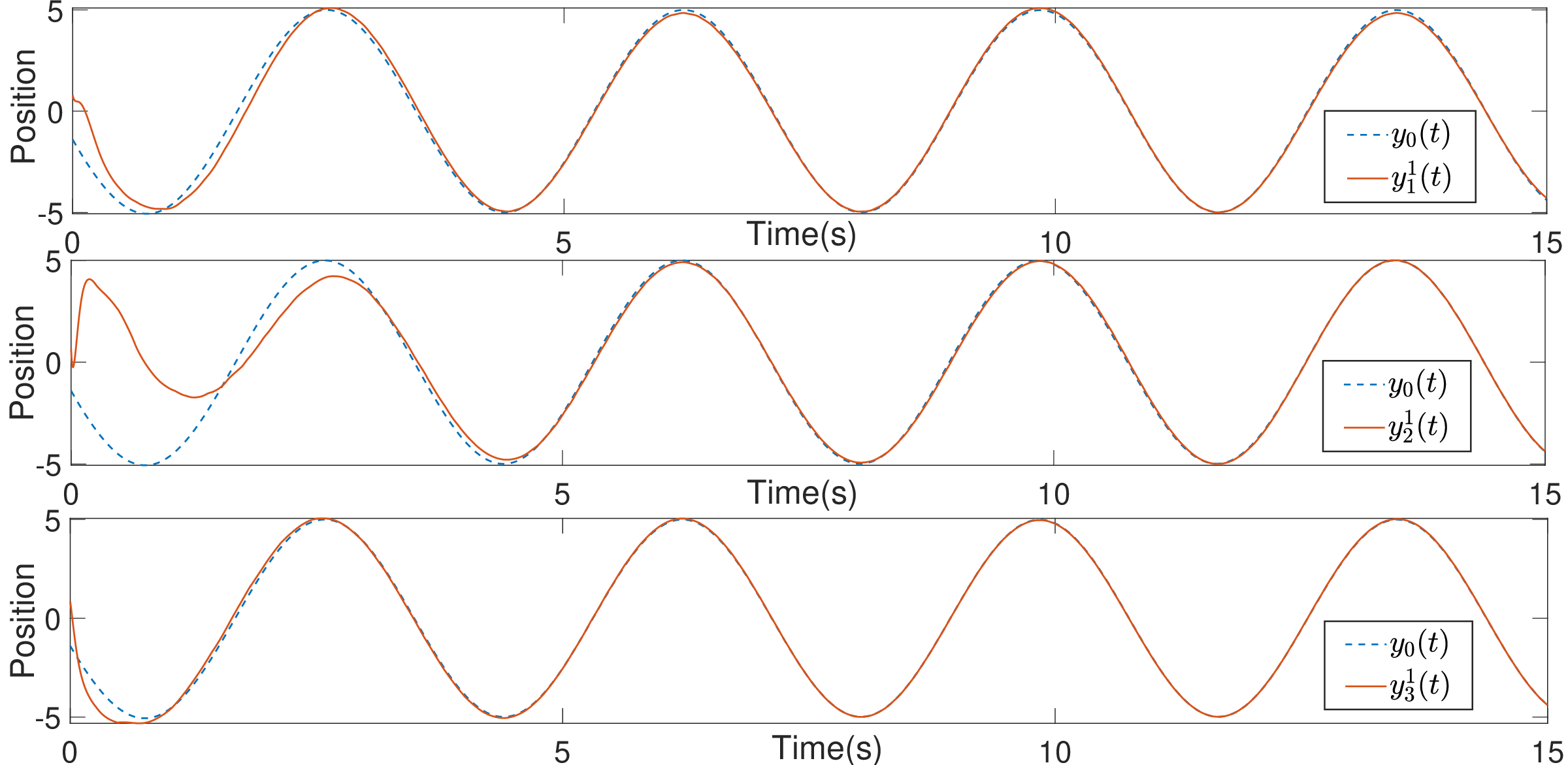}
\caption{The trajectory tracking of the leader and the $1^{st}$ type followers.}
\label{Pi_1_y1}
\end{figure}

\begin{figure}[!ht]
\centering
\includegraphics[width=0.7\linewidth,height = 4.5cm]{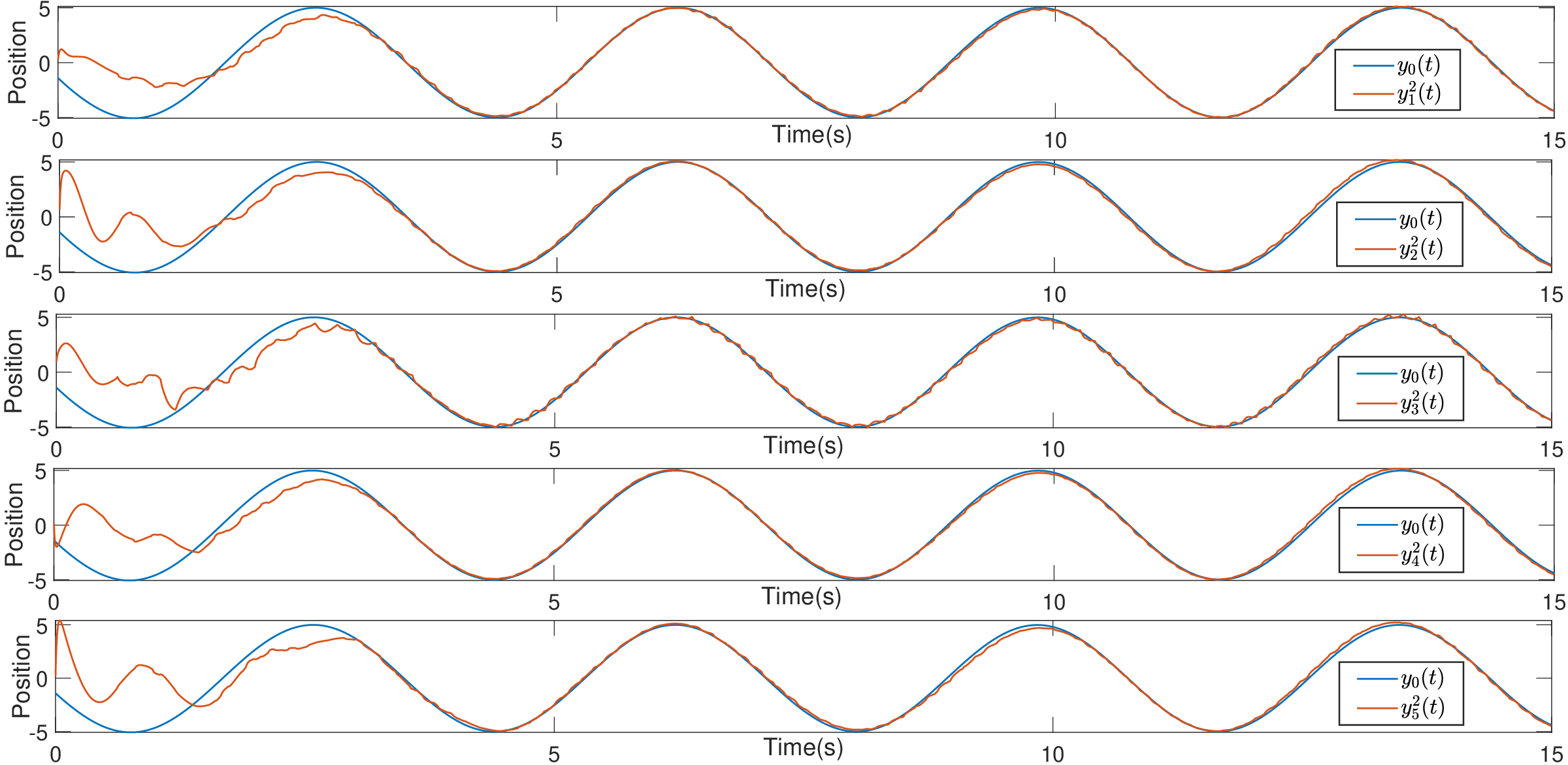}
\caption{The trajectory tracking of the leader and the $2^{nd}$ type followers.}
\label{Pi_2_y_wave}
\end{figure}

The trajectories of the observer errors $||\tilde x_i^1|{|^2}$ and $||\tilde x_j^2|{|^2}$ are demonstrated in Fig.\ref{Pi_1_x_wave} and Fig.\ref{Pi_2_x_wave}, respectively. The errors $||\varepsilon _i^1|{|^2}$ and $||\varepsilon _j^2|{|^2}$ are illustrated in Fig.\ref{Pi_1_eplison} and Fig.\ref{Pi_2_eplison} respectively, $i \in {\Lambda _1}$, $j \in {\Lambda _2}$. The event-triggered intervals (ETI) of different agents are shown in Figs.\ref{Pi1_event}-\ref{Pi2_event}. Then the estimation of ${A_r}$ is displayed in Fig.\ref{Pi_2_A0}. Based on the distributed controllers in Eqs.(\ref{type1_11}) and (\ref{type2_20}), the trajectory tracking curves of the followers and the leader are drawn in Fig.\ref{Pi_1_y1} and Fig.\ref{Pi_2_y_wave}. It can be seen that the leader-follower consensus can be achieved with UUB via the proposed control strategy. By comparing the state estimation and trajectory tracking of the two types followers, it can be observed that there is no significant difference between the two types of followers. 
Hence, it can demonstrate that the proposed adaptive event-triggered mechanism can be used to achieve control consistency in the general heterogeneous multi-agent system.

\section{Conclusion}
In this paper, the consistency problem of the heterogeneous multi-agent system is studied with simultaneous communication faults and actuator faults. The followers are divided into two groups, directly and not directly connected with the leader, where the leader model can be influenced by the $1^{st}$ type followers. Since the $2^{nd}$ type followers can only communicate with the $1^{st}$ type followers, the system matrix estimators are designed for the $2^{nd}$ type followers to estimate the matrix of the leader. A series of new distributed event-triggered observers are designed for the $2^{nd}$ types followers to estimate the leader state. In order to reduce the communication burden, an event-triggered controller is proposed for each agent, which can ensure that the leader-followers is consensus with UUB for the studied heterogeneous multi-agent system. Simulation experiments have been performed to verify the effectiveness of the proposed method. When the studied system approaches infinity, it can achieve consensus. However, this challenges the practical applications. Therefore, the consensus issue will be further explored with limited time to accommodate in more practical scenarios. Therefore, the next step is to explore the issue of consensus in systems with predefined time. The consensus issue will be further explored with prescribed time to accommodate in more practical scenarios.

%\small{
%\begin{ack}                               % Place acknowledgements
%This work was supported in part by the National Natural Science Foundation of China Ref. U1913201 and Ref. 61973296, in part by the Guangdong Basic and Applied Basic Research Foundation Ref.2021B1515120038, in part by the Shenzhen Science and Technology Innovation Commission Project Grant Ref.JCYJ20200109114839874, Ref.JSGG20210802154535003 and Ref.JCYJ20210324101215039.  % here.
%\end{ack}
%}

\appendix
\section{Proof of Lemma 1}
Different gain matrixes $W_i^1$ are designed for the $1^{st}$ type followers, so that $A_i^1 - W_i^1$ are Hurwitz matrices. According to Definition 1, there is a positive definite matrix ${Q_i}$, $i = 1, \cdots m$, so that
$${(A_i^1 - W_i^1)^T}{Q_i} + {Q_i}(A_i^1 - W_i^1) <  - {\alpha _i}{Q_i},$$
where ${\alpha _i}$ is a positive real value, $i = 1, \cdots m$.
Considering a Lyapunov function,
\begin{equation} \label{leader_nof6}
{V_1} = {\bar x}^{1T}Q{\bar x^1} + {\tilde u}^{{a_1T}}{\tilde u}^{{a_1}},
\end{equation}
where $Q = [{Q_1}, \cdots ,{Q_m}]$, ${V_1} = diag({V_{11}}, \cdots ,{V_{1m}})$, ${V_{1i}} = \bar x_i^{1T}{Q_i}\bar x_i^1 + \tilde u{_i^{a_1 T}}\tilde u_i^{{a_1}}$.

Taking the time derivative of ${V_1}$ along the trajectories of System (\ref{leader_nof5}),
\begin{equation} \label{type1_7}
\begin{aligned}
{{\dot V}_1} =& {{\bar x}^{1T}}({({A^1} - {W^1})^T}Q + Q({A^1} - {W^1})){{\bar x}^1}+ 2{{\tilde u}^{{a_1}T}}{{\dot u}^{{a_1}}}\\
&+2{{\bar x}^{1T}}Q^T{B^{1}}{{\tilde u}^{{a_1}}}+ 2{{\tilde u}^{{a_1}T}}{N^{11}}{{\hat u}^{{a_1}}} - 2{{\tilde u}^{{a_1}T}}{N^{12}}(x^1 - \hat x^1)\\
\le&  -{{\bar x}^{1T}}{\alpha ^1}Q{{\bar x}^1} + 2{{\tilde u}^{{a_1}T}}{B^1}^TQ{{\bar x}^1}
+ 2{{\tilde u}^{{a_1}T}}{{\dot u}^{{a_1}}} \\
&+ 2{{\tilde u}^{{a_1}T}}{N^{11}}{{\hat u}^{{a_1}}}
- 2\tilde u{_i^{{a_1}T}}{N^{12}}(x^1 - \hat x^1).
\end{aligned}
\end{equation}
where ${\alpha ^1}= ({\hat \alpha}^1 \otimes I_n), {{\hat \alpha}^1} = diag({\alpha_{1}}, \cdots ,{\alpha_{m}})$.

According to Eq.(\ref{leader_nof5.1}), it has,
$$2{\tilde u^{{a_1}T}}{N^{11}}{\hat u^{{a_1}}} = 2{\tilde u^{{a_1}T}}{N^{11}}{u^{{a_1}}} - 2{\tilde u^{{a_1}T}}{N^{11}}{\tilde u^{{a_1}}}.$$

By Young's inequality,
\begin{align}
\label{type1_8} 2{{\tilde u}^{{a_1}T}}{{\dot u}^{{a_1}}} & \le{{\tilde u}^{{a_1}T}}{\alpha ^2}{{\tilde u}^{{a_1}}} +{{\dot u}^{{a_1}T}}\frac{1}{{{\alpha ^2}}}{{\dot u}^{{a_1}}},\\
\label{type1_9} 2{{\tilde u}^{{a_1}T}}{N^{11}}{u^{{a_1}}} &\le{{\tilde u}^{{a_1}T}}{\alpha ^3}{N^{11}}{{\tilde u}^{{a_1}}} +{u^{{a_1}T}}\frac{1}{{{\alpha ^3}}}{N^{11}}{u^{{a_1}}},
\end{align}
where ${\alpha ^2}$ and ${\alpha ^3}$ are the constants to be determined. Let ${N^{12}} = {Q^T}{B^1}$.

Substituting Eq.(\ref{type1_8}) and Eq.(\ref{type1_9}) into Eq.(\ref{type1_7}), it has,
\begin{equation}\label{type1_10}
\begin{aligned}
{{\dot V}_1} \le&  -{{\bar x}^{1T}}{\alpha ^1}Q{{\bar x}^1} + {{\tilde u}^{{a_1}T}}{\alpha ^2}{{\tilde u}^{{a_1}}}+{{\tilde u}^{{a_1}T}}{\alpha ^3}{N^{11}}{{\tilde u}^{{a_1}}} \\
& +{{\dot u}^{{a_1}T}}\frac{1}{{{\alpha ^2}}}{{\dot u}^{{a_1}}}
+{u^{{a_1}T}}\frac{1}{{{\alpha ^3}}}{N^{11}}{u^{{a_1}}} - 2{{\tilde u}^{{a_1}T}}{N^{11}}{{\tilde u}^{{a_1}}}.
\end{aligned}
\end{equation}

With appropriate ${\alpha ^2}$, ${\alpha ^3}$ and ${N^{11}}$ such that
$${\alpha ^2} + {\alpha ^3}{N^{11}} - {N^{11}} < 0.$$

Then
\begin{equation}
\begin{aligned}
{\dot V_1} \le & - {\bar x^{1T}}{\alpha ^1}Q{\bar x^1} +{\dot u^{{a_1}T}}\frac{1}{{{\alpha ^2}}}{\dot u^{{a_1}}}
+{u^{{a_1}T}}\frac{1}{{{\alpha ^3}}}{N^{11}}{u^{{a_1}}} - {\tilde u^{{a_1}T}}{N^{11}}{\tilde u^{{a_1}}}.
\end{aligned}
\end{equation}

According to Assumption 4, ${u^{{a_1}}}$ and ${\dot u^{{a_1}}}$ is bounded, then
\begin{equation}
\begin{aligned}
{{\dot V}_1}& \le  - {{\bar x}^{1T}}{\alpha ^1}Q{{\bar x}^1} + {\varphi _1} - {{\tilde u}^{{a_1}T}}{N^{11}}{{\tilde u}^{{a_1}}}\\
& \le  - {\varphi _2}({{\bar x}^{1T}}Q{{\bar x}^1} - {{\tilde u}^{{a_1}T}}{{\tilde u}^{{a_1}}}) + {\varphi _1}\\
& = - {\varphi _2}{V_1} + {\varphi _1},
\end{aligned}
\end{equation}
where ${\varphi _1} = \max \{{\dot u^{{a_1}T}}\frac{1}{{{\alpha ^2}}}{\dot u^{{a_1}}} +{u^{{a_1 T}}}\frac{1}{{{\alpha ^3}}}{N^{11}}{u^{{a_1}}}\}$, ${\varphi _2} = \min \{ \alpha _1^1, \cdots ,\alpha _m^1,{\lambda _{\min }}({N^{11}})\}$.
Then,
\begin{equation}
\begin{array}{l}
{V_{1i}} < ({V_{1i}}(0) - \frac{{{\varphi _1}}}{{{\varphi _2}}}){e^{ - {\varphi _2}t}} + \frac{{{\varphi _1}}}{{{\varphi _2}}},
\end{array}
\end{equation}
when $t \to \infty$, ${e^{ - {\varphi _2}t}}$ approaches 0, indicating that if $({V_i}(0) - \frac{{{\varphi _1}}}{{{\varphi _2}}})$ is non-zero, for $\forall t \ge {T_i}$, there is ${e^{ - {\varphi _2}t}} < \frac{{{\varphi _1}}}{{{\varphi _2}({V_i}(0) - \frac{{{\varphi _1}}}{{{\varphi _2}}})}}$. According to Eq.(\ref{leader_nof6}), for $\forall t \ge {T_i}$, it is known that $||\bar x_i^1|| < \lambda _{\min }^{ - 1}({Q_i}){V_{1i}}$, $||u_i^{{a_1}}|| < {V_{1i}}$, $i = 1, \cdots ,m$. Thus $\bar x_i^1$, $\tilde u_i^{a_1}$ are UUB.   $\hfill{} \Box$

\section{Proof of Lemma 7}
Let
$\bar v = {[\bar v_{m + 1}^T, \cdots ,\bar v_N^T]^T} =  - {\varphi ^2}({H_2} \otimes {I_{{n^2}}})(\mathord{\buildrel{\lower3pt\hbox{$\scriptscriptstyle\smile$}}
\over A}  + \bar e)$, it has,
\begin{equation}\label{type2_5}
\begin{array}{l}
\mathord{\buildrel{\lower3pt\hbox{$\scriptscriptstyle\smile$}}
\over A}  =  - \frac{1}{{{\varphi ^2}}}(H{}_2^{ - 1} \otimes {I_{{n^2}}})\bar v - \bar e.
\end{array}
\end{equation}

Considering the Lyapunov function,
\begin{equation}\label{type2_5._1}
\begin{array}{l}
{V_3} = {\mathord{\buildrel{\lower3pt\hbox{$\scriptscriptstyle\smile$}}
\over {A^T}}}({H_2} \otimes {I_{{n^2}}})\mathord{\buildrel{\lower3pt\hbox{$\scriptscriptstyle\smile$}}
\over A}  + \sum\limits_{i = m + 1}^N {{\phi _{i2}}}.
\end{array}
\end{equation}

The derivative of ${V_3}$ along the trajectories of Eq.(\ref{type2_4}) is,
\begin{equation}\label{type2_5.0}
\begin{array}{l}
{\dot V_3} = 2{\mathord{\buildrel{\lower3pt\hbox{$\scriptscriptstyle\smile$}}
\over A} ^T}({H_2} \otimes {I_{{n^2}}})\bar v + {{\tilde A}^T}(t)(\dot H \otimes {I_{{n^2}}})\tilde A(t) + \sum\limits_{i = m + 1}^N {{{\dot \phi }_{i2}}}.
\end{array}
\end{equation}

From Eq.(\ref{type2_5}), it has,
\begin{equation} \label{type2_5.1}
\begin{aligned}
{{\dot V}_3} \le&  - (\frac{2}{{\varphi _{\max }^2}} - {\lambda _1})\sum\limits_{j = m + 1}^N {\bar v_j^T{{\bar v}_j}}  - 2\sum\limits_{j = m + 1}^N {d_j^2\bar e_j^T{{\bar v}_j}} 
 + 2\sum\limits_{i = m + 1}^N {\sum\limits_{j = m + 1}^N {{{\bar a}_{ij}}\bar e_i^T{{\bar v}_j}} } \\
& - 2\sum\limits_{j = m + 1}^N {d_j^{21}\bar e_j^T{{\bar v}_j}}
+ \sum\limits_{i = m + 1}^N {{{\dot \phi }_{i2}}}  + 2{\lambda _2}{{\bar e}^T}\bar v + {\lambda _3}\sum\limits_{j = m + 1}^N {\bar e_j^T{{\bar e}_j}}, 
\end{aligned}
\end{equation}
where ${\lambda _1} = {\lambda _{\max }}(H_2^{ - 1}\frac{1}{{{\varphi ^2}}}{\dot H_2}\frac{1}{{{\varphi ^2}}}H_2^{ - 1} \otimes {I_{{n^2}}})$, ${\lambda _2} = 2{\lambda _{\max }}({\dot H_2}\frac{1}{{{\varphi ^2}}}H_2^{ - 1} \otimes {I_{{n^2}}})$, ${\lambda _3} = {\lambda _{\max }}({\dot H_2} \otimes {I_{{n^2}}})$, $\varphi _{\max }^2 = \max \{ \varphi _{m + 1}^2, \cdots ,\varphi _N^2\}$.

By Young's inequality,
\begin{align}
\label{type2_6} 2\sum\limits_{i = m + 1}^N {\sum\limits_{j = m + 1}^N {{{\bar a}_{ij}}\bar e_i^T{{\bar v}_j}} } &\le \sum\limits_{i = m + 1}^N {\sum\limits_{j = m + 1}^N {{{\bar a}_{ij}}({{\bar k}_{1i}}\bar e_i^T{{\bar e}_i} + \frac{1}{{{{\bar k}_{1i}}}}\bar v_j^T{{\bar v}_j}} }),\\
\label{type2_7} - 2\sum\limits_{j = m + 1}^N {d_j^2\bar e_j^T{{\bar v}_j}}  &\le 
\sum\limits_{j = m + 1}^N {d_j^2({{\bar k}_{1j}}\bar e_j^T{{\bar e}_j} + \frac{1}{{{{\bar k}_{1j}}}}\bar v_j^T{{\bar v}_j}} ),\\
\label{type2_8} - 2\sum\limits_{j = m + 1}^N {d_j^{21}\bar e_j^T{{\bar v}_j}} &\le 
\sum\limits_{j = m + 1}^N {d_j^{21}({{\bar k}_{1j}}\bar e_j^T{{\bar e}_j} + \frac{1}{{{{\bar k}_{1j}}}}\bar v_j^T{{\bar v}_j})},\\
\label{type2_8_1}2{\lambda _2}\sum\limits_{j = m + 1}^N {\bar e_j^T{{\bar v}_j}} &\le {\lambda _2}\sum\limits_{j = m + 1}^N {({{\bar k}_{1j}}\bar e_j^T{{\bar e}_j} + \frac{1}{{{{\bar k}_{1j}}}}\bar v_j^T{{\bar v}_j})},
\end{align}
where ${\bar k_{1j}} > 0$. Since the communications among the followers of the $2^{nd}$ group are undirected, it has,
\begin{equation} \label{type2_9}
\begin{aligned}
\sum\limits_{i = m + 1}^N {\sum\limits_{j = m + 1}^N {{{\bar a}_{ij}}{{\bar k}_{1j}}\bar e_i^T{{\bar e}_i}} } &= \sum\limits_{i = m + 1}^N {\sum\limits_{j = m + 1}^N {{{\bar a}_{ij}}{{\bar k}_{1i}}\bar e_i^T{{\bar e}_i}}} 
\le \sum\limits_{i = m + 1}^N {(d_i^2 + d_i^{21}){{\bar k}_{1i}}\bar e_i^T{{\bar e}_i}},
\end{aligned}
\end{equation}

Taking Eqs.(\ref{type2_6})-(\ref{type2_9}) into Eq.(\ref{type2_5.1}), it has,
%\begin{equation} \label{type2_9_1}
%\begin{array}{l}
%{\dot V_3} \le \sum\limits_{j = m + 1}^N {\left[ { - {{\hat \alpha }_j}{{\left\| {{{\bar v}_j}} \right\|}^2} + 2\left( {d_j^2 + d_j^{21}} \right){{\bar k}_{1j}}{{\left\| {{{\bar e}_j}} \right\|}^2}} \right]} \\
% + \sum\limits_{i = m + 1}^N {{{\dot \phi }_{i2}}}
%\end{array}
%\end{equation}

\begin{equation} \label{type2_9_1}
\begin{array}{l}
{\dot V_3} \le \sum\limits_{j = m + 1}^N {\left[ { - {{\hat \alpha }_j}{{\left\| {{{\bar v}_j}} \right\|}^2} + {{\hat \beta }_j}{{\left\| {{{\bar e}_j}} \right\|}^2}} \right]}  + \sum\limits_{i = m + 1}^N {{{\dot \phi }_{i2}}},
\end{array}
\end{equation}
where ${\hat \alpha _j} =  - ( - \frac{2}{{\varphi _{\max }^2}} + {\lambda _1} + \frac{{d_j^2}}{{{{\bar k}_{1j}}}} + \frac{{d_j^{21}}}{{{{\bar k}_{1j}}}} + \frac{{{\lambda _2}}}{{{{\bar k}_{1j}}}} + \sum\limits_{i = m + 1}^N {\frac{{{{\bar a}_{ij}}}}{{{{\bar k}_{1i}}}}} )$, ${\hat \beta _j} = (d_j^2 + d_j^{21} + {\lambda _2} + \frac{{{\lambda _3}}}{{{{\bar k}_{1j}}}}){\bar k_{1j}}$. Then the appropriate $\varphi _i^2$ and ${\bar k_{1j}}$ parameters can be set so that ${\hat \alpha _j} \ge 0$, $i,j = m + 1, \cdots ,N$. From Eq.(\ref{type2_3}), it has,
\begin{equation} \label{type2_9_2}
\begin{array}{l}
{\left\| {{{\bar e}_i}} \right\|^2} \le \frac{1}{{{{\hat \omega }_i}{\gamma _{i2}}}}{\phi _{i2}} - \frac{1}{{{\theta _{i2}}{\gamma _{i2}}}}||{\psi _{i2}}|{|^2} + \frac{\hat\delta _{i2}}{\gamma _{i2}}.
\end{array}
\end{equation}
where ${\hat\delta _{i2}} = {\tau _{i2}}{e^{ - \frac{{{\delta _{i2}}}}{{{\sigma _{i2}} + t}}}} + {{\hat \tau }_{i2}}$.
Then,
\begin{small}
\begin{equation} \label{type2_10}
\begin{aligned}
{{\dot V}_3} \le& \sum\limits_{j = m + 1}^N {\left[ {{{\hat \beta }_j}\left( {\frac{1}{{{{\hat \omega }_j}{\gamma _{j2}}}}{\phi _{j2}} - \frac{1}{{{\theta _{j2}}{\gamma _{j2}}}}||{\psi _{j2}}|{|^2} + \frac{\hat\delta _{j2}}{\gamma _{j2}}} \right)} \right]} \\
 &+ \sum\limits_{i = m + 1}^N {\left( { - {\beta _{i2}}{\phi _{i2}} - {\gamma _{i2}}||{{\bar e}_i}|{|^2} + \frac{1}{{{\theta _{i2}}}}||{\psi _{i2}}|{|^2} + {\hat\delta _{i2}}} \right)} \\
 \le &\sum\limits_{j = m + 1}^N {\hat \alpha _j^1{\phi _{j2}}}  - \hat \alpha _{\min }^2\sum\limits_{j = m + 1}^N {||{\psi _{j2}}|{|^2}}  + \sum\limits_{j = m + 1}^N {\hat \alpha _j^3}, 
\end{aligned}
\end{equation}
\end{small}
where $\hat \alpha _j^1 = \frac{{{{\hat \beta }_j}}}{{{{\hat \omega }_j}{\gamma _{j2}}}} - {\beta _{j2}}$, $\hat \alpha _{j}^2=\frac{{{{\hat \beta }_j}}}{{{\theta _{j2}}{\gamma _{j2}}}} - \frac{1}{{{\theta _{j2}}}}$, $\hat \alpha _j^3 = \frac{{{\hat \beta }_j}{\hat\delta _{j2}}}{\gamma _{j2}} + {\hat\delta _{j2}}$ and $\hat \alpha _{{\rm{min}}}^2 = \max \{ \hat \alpha _{m + 1}^2, \cdots ,\hat \alpha _N^2\}$. Then,
\begin{equation} \label{type2_11}
\begin{array}{l}
- \sum\limits_{j = m + 1}^N {||{\psi _{j2}}|{|^2} \le  - ||({H_2} \otimes {I_{{n^2}}})({{\bar A}_r} - \bar A)|{|^2}} 
\le  - \frac{1}{{\varphi _{\max }^2}}||\bar v|{|^2}.
\end{array}
\end{equation}

Taking Eq.(\ref{type2_11}) into Eq.(\ref{type2_10}), it has,
\begin{equation} \label{type2_11_1}
\begin{array}{l}
{\dot V_3} \le \sum\limits_{j = m + 1}^N {\hat \alpha _j^1{\phi _{j2}}}  - \frac{{\hat \alpha _{\min }^2}}{{\varphi _{\max }^2}}||\bar v|{|^2} + \sum\limits_{j = m + 1}^N {\hat \alpha _j^3}.
\end{array}
\end{equation}

According to Assumption 1 and Assumption 2, ${H_2}$ is a positive definite matrix. From Eq.(\ref{type2_5}) and Eq.(\ref{type2_5._1}), it has,
\begin{equation}\label{type2_12}
\begin{aligned}
{V_3} =& {( - \frac{1}{{{\varphi ^2}}}(H{}_2^{ - 1} \otimes {I_{{n^2}}})\bar v - \bar e)^T}({H_2} \otimes {I_{{n^2}}})
( - \frac{1}{{{\varphi ^2}}}(H{}_2^{ - 1} \otimes {I_{{n^2}}})\bar v - \bar e) + \sum\limits_{i = m + 1}^N {{\phi _{i2}}} \\
\le & \frac{{{\lambda _{\max }}(H)}}{{\varphi _{\min }^2\varphi _{\min }^2\lambda _{\min }^2(H)}}||\bar v|{|^2} + \bar e{^T}({H_2} \otimes {I_{{n^2}}})\bar e 
 + 2\bar e{^T}({H_2} \otimes {I_{{n^2}}})\frac{1}{{{\varphi ^2}}}(H{}_2^{ - 1} \otimes {I_{{n^2}}})\bar v 
+ \sum\limits_{i = m + 1}^N {{\phi _{i2}}}, 
\end{aligned}
\end{equation}
where $\varphi _{\min }^2 = \min \{ \varphi _{m + 1}^2, \cdots ,\varphi _N^2\}$.
By Young's inequality,
\begin{equation}\label{type2_13}
\begin{array}{l}
2\bar e{^T}({H_2} \otimes {I_{{n^2}}})\frac{1}{{{\varphi ^2}}}(H{}_2^{ - 1} \otimes {I_{{n^2}}})\bar v 
\le \frac{{{{\tilde k}_1}}}{{\varphi _{\min }^2}}||\bar v|{|^2} + \frac{1}{{\varphi _{\min }^2{{\tilde k}_1}}}||\bar e|{|^2},
\end{array}
\end{equation}
where ${\tilde k_1} > 0$. Taking Eq.(\ref{type2_13}) into Eq.(\ref{type2_12}), and from Eq.(\ref{type2_9}), it has,
\begin{equation}\label{type2_14}
\begin{aligned}
{V_3} \le &\frac{{{\lambda _{\max }}(H)}}{{\varphi _{\min }^2\varphi _{\min }^2\lambda _{\min }^2(H)}}||\bar v|{|^2} + \bar e{^T}({H_2} \otimes {I_{{n^2}}})\bar e 
 + 2\bar e{^T}({H_2} \otimes {I_{{n^2}}})\frac{1}{{{\varphi ^2}}}(H{}_2^{ - 1} \otimes {I_{{n^2}}})\bar v + \sum\limits_{i = m + 1}^N {{\phi _{i2}}} \\
\le & {{\hat \gamma }_2}\sum\limits_{i = m + 1}^N {\left( {\frac{1}{{{{\hat \omega }_i}{\gamma _{i2}}}}{\phi _{i2}} - \frac{1}{{{\theta _{i2}}{\gamma _{i2}}}}||{\psi _{i2}}|{|^2}} + \frac{{\hat \alpha _i^4} +\phi _{i2}} {{\hat\gamma_{2}}}\right)}+{{\hat \gamma }_1}||\bar v|{|^2}\\
\le& \left( {{{\hat \gamma }_1} - \frac{{{{\hat \gamma }_2}}}{{\varphi _{\max }^2{{\hat \gamma }_3}}}} \right)||\bar v|{|^2} + \sum\limits_{i = m + 1}^N \left({\frac{{{\phi _{i2}}{{\hat \gamma }_2}}}{{{{\hat \omega }_i}{\gamma _{i2}}}} + {\phi _{i2}}} + {\hat \alpha _i^4}\right), \\
\end{aligned}
\end{equation}
where ${\hat \gamma _1} = \frac{{{\lambda _{\max }}({H_2})}}{{\varphi _{\min }^2\varphi _{\min }^2\lambda _{\min }^2({H_2})}} + \frac{{{{\tilde k}_1}}}{{\varphi _{\min }^2}}$, ${\hat \gamma _2} = \frac{1}{{\varphi _{\min }^2{{\tilde k}_1}}} + {\lambda _{\max }}({H_2})$, ${\hat \gamma _3} = \max \{ {\theta _{(m + 1)2}}{\gamma _{(m + 1)2}}, \cdots ,{\theta _{N2}}{\gamma _{N2}}\}$ and ${\hat \alpha _i^4}=\frac{{{\hat \gamma }_2}{\hat \delta _{i2}}}{{{\gamma _{i2}}}}$. 

Let ${\hat \gamma _4} = {\hat \gamma _1} - \frac{{{{\hat \gamma }_2}}}{{\varphi _{\max }^2{{\hat \gamma }_3}}} > 0$.

From Eq.(\ref{type2_14}), it has,
\begin{equation}\label{type2_15}
\begin{array}{l}
\frac{1}{{{{\hat \gamma }_4}}}{V_3} - \frac{1}{{{{\hat \gamma }_4}}}\sum\limits_{i = m + 1}^N\left( {\frac{{{\phi _{i2}}{{\hat \gamma }_2}}}{{{{\hat \omega }_i}{\gamma _{i2}}}} + {\phi _{i2}}} + {\hat \alpha _i^4}\right) \le ||\bar v|{|^2}.
\end{array}
\end{equation}

From Eq.(\ref{type2_11_1}) and Eq.(\ref{type2_15}), it has,
\begin{equation}\label{type2_16}
\begin{aligned}
{\dot V_3} \le& \sum\limits_{j = m + 1}^N {\hat \alpha _j^1{\phi _{j2}}}  - \frac{{\hat \alpha _{\min }^2}}{{{{\hat \gamma }_4}}{\varphi _{\max }^2}}{V_3} + \sum\limits_{j = m + 1}^N {\hat \alpha _j^3}
+\frac{{\hat \alpha _{\min }^2}}{{{{\hat \gamma }_4}}{\varphi _{\max }^2}}\sum\limits_{i = m + 1}^N\left( {\frac{{{\phi _{i2}}{{\hat \gamma }_2}}}{{{{\hat \omega }_i}{\gamma _{i2}}}} + {\phi _{i2}}} + {\hat \alpha _i^4}\right)\\
=&  - {\eta ^1}{V_3} + \sum\limits_{i = m + 1}^N {\eta _i^2} {\phi _{i2}} + {\eta^3}.
\end{aligned}
\end{equation}
where ${\eta ^1} = \frac{{\hat \alpha _{\min }^3}}{{{\hat \gamma }_4}{\varphi _{\max }^2}}$, $\eta _i^2 = \hat \alpha _i^1 + \frac{{\hat \alpha _{\min }^2}}{{{\hat \gamma }_4}{\varphi _{\max }^2}}\left( {\frac{{{{\hat \gamma }_2}}}{{{{\hat \omega }_i}{\gamma _{i2}}}} + 1} \right)$ and $\eta^3 = \sum\limits_{i = m + 1}^N \left(\frac{{\hat \alpha _{\min }^2}{\hat \alpha _{i}^4}}{{{\hat \gamma }_4}{\varphi _{\max }^2}} + {\hat \alpha _i^3}\right)$,

Then appropriate parameters ${\hat \omega_j}$, ${\gamma_{j2}}$ and ${\beta_{j2}}$ can be set so that $\eta_i^2 \leq 0$, $i = m + 1, \cdots ,N$.
Then it has,
\begin{equation}\label{type2_17}
\begin{array}{l}
{\dot V_3} \le - {\eta ^1}{V_3} + {\eta^3}.
\end{array}
\end{equation}

It can be obtained that ${V_3} \le ({V_3}(0)-\frac{\eta^3}{\eta^1}){e^{ - {\eta ^1}t}}+\frac{\eta^3}{\eta^1}$, that is, when $t \to \infty$, ${e^{ - {\eta ^1}t}}$ approaches 0. According to Eq.(\ref{type2_5._1}), it is known that ${\lambda _{\min }}({H_2})||\mathord{\buildrel{\lower3pt\hbox{$\scriptscriptstyle\smile$}}
\over A} |{|^2} \le {V_3}$, which indicates that $\mathord{\buildrel{\lower3pt\hbox{$\scriptscriptstyle\smile$}}
\over A} $ is bounded. Therefore, it can be concluded that ${\mathord{\buildrel{\lower3pt\hbox{$\scriptscriptstyle\smile$}}
\over A} _i} = ({A_r} - {\hat A_i})$  reach the consistency with UUB. $\hfill{} \Box$

\section{Proof of Lemma 11}

Considering the following Lyapunov function,
$${W_2} = {V_2} + {\hat V_2}.$$

From Eq.(\ref{type1_27}) and Eq.(\ref{type2_26}), it has,
\begin{small}
\begin{equation}
\begin{aligned}
{{\dot W}_2} = & {{\dot V}_2} + {{\dot {\hat V}}_2} = - {{\hat \alpha }_2}{c_1}{{\tilde x}^{1T}}{P_1}{{\tilde x}^1} + \frac{1}{{{\beta _2}}}{{\tilde x}^{2T}}{{\tilde x}^2}\\
& + \frac{1}{{{\tau _2}}}||{P_2}X{M^{11}}{M^{13}}({A^{12}} \otimes {I_n}){{\tilde x}^2}|{|^2} - {{\tilde \alpha }_1}{\varepsilon ^1}^T{_1}{P_2}{\varepsilon ^1}\\
& + \frac{1}{{{\beta _4}}}||{{\bar g}_1}{{\hat K}_1}{\varepsilon ^1}|{|^2} + {\tau _6}{\varepsilon ^1}^T{_1}{\varepsilon ^1} + {\Omega _1} + {\Omega _3} + {{\hat \Omega }_1} + {{\hat \Omega }_2}\\
& + \hat M - \sum\limits_{j = 1}^m ({{\beta _{j1}}} {\phi _{j1}}-{\hat \delta_{j1}}) - {\mu _2}{c_2}{{\tilde x}^{2T}}{{\hat P}_1}{{\tilde x}^2}\\
& + {\beta _3}||{{\tilde x}^{1T}}{K^1}|{|^2} + \frac{1}{{{{\hat \beta }_2}}}||{M^{21}}{M^{23}}({A^{21}} \otimes {I_n}){{\tilde x}^1}|{|^2}\\
& + {\beta _4}||{{\tilde x}^{1T}}{K^1}|{|^2} - {{\tilde \mu }_1}{\varepsilon ^2}^T{{\hat P}_2}{\varepsilon ^2} + {{\hat \tau }_6}{\varepsilon ^2}^T{\varepsilon ^2}\\
& + \tilde M + \frac{1}{{{{\hat \tau }_2}}}||{{\hat P}_2}\hat X{M^{21}}{M^{23}}({A^{21}} \otimes {I_{(N - m)}}){{\tilde x}^1}|{|^2}\\
& + \frac{1}{{{{\hat \tau }_3}}}||{{\hat P}_2}\hat X{M^{21}}{M^{23}}({A^{21}} \otimes {I_{(N - m)}}){{\tilde x}^2}|{|^2}\\
&- \sum\limits_{i = m + 1}^N {\left( {{\beta _{i3}}{\phi _{i3}} -{\hat \delta_{j3}}} \right)}.
\end{aligned}
\end{equation}
\end{small}

Let set appropriate parameters ${\hat \alpha _1}$, ${K^1}$, ${M^{11}}$, ${M^{13}}$, ${\beta _1}$, ${\beta _3}$, ${\beta _4}$, ${M^{21}}$, ${M^{23}}$, ${\hat \beta _1}$, ${\hat \beta _3}$ and ${\hat \beta _4}$, so that
\begin{small}
\begin{equation} \label{type2_27}
\begin{array}{l}
\frac{1}{{{{\hat \beta }_2}}}||{M^{21}}{M^{23}}({A^{21}} \otimes {I_n}){{\tilde x}^1}|{|^2} + {\beta _4}||{{\tilde x}^{1T}}{K^1}|{|^2} + {\Omega _3} + \\
{\beta _3}||{{\tilde x}^{1T}}{K^1}|{|^2} + \frac{1}{{{{\hat \tau }_2}}}||{{\hat P}_2}\hat X{M^{21}}{M^{23}}({A^{21}} \otimes {I_{(N - m)}}){{\tilde x}^1}|{|^2} < 0,
\end{array}
\end{equation}
\begin{equation} \label{type2_28}
\begin{array}{l}
- {\hat c_1}{\tilde \alpha _1}{\varepsilon ^1}^T{_1}{P_2}{\varepsilon ^1} + \frac{1}{{{\beta _4}}}||{\bar g_1}{\hat K_1}{\varepsilon ^1}|{|^2} + {\tau _6}{\varepsilon ^1}^T{_1}{\varepsilon ^1} < 0,
\end{array}
\end{equation}
\begin{equation} \label{type2_29}
\begin{array}{l}
{\hat \Omega _2} + \frac{1}{{{\beta _2}}}{\tilde x^{2T}}{\tilde x^2} + \frac{1}{{{\tau _2}}}||{P_2}X{M^{11}}{M^{13}}({A^{12}} \otimes {I_n}){\tilde x^2}|{|^2} \\
+ \frac{1}{{{{\hat \tau }_3}}}||{\hat P_2}\hat X{M^{21}}{M^{23}}({A^{21}} \otimes {I_{(N - m)}}){\tilde x^2}|{|^2} < 0,
\end{array}
\end{equation}
\begin{equation} \label{type2_30}
\begin{array}{l}
- {\hat c_2}{\tilde \mu _1}{\varepsilon ^2}^T{\hat P_2}{\varepsilon ^2} + {\hat \tau _6}{\varepsilon ^2}^T{\varepsilon ^2} < 0,
\end{array}
\end{equation}
\end{small}
where $0 < {\hat c_1},{\hat c_2} < 1$.

From Eqs.(\ref{type2_27})-(\ref{type2_30}), it has,
\begin{equation}\label{type2_31}
\begin{aligned}
{{\dot W}_2} \le & - {{\hat \alpha }_2}{c_1}{{\tilde x}^{1T}}{P_1}{{\tilde x}^1} - {{\tilde \alpha }_1}{\varepsilon ^1}^T{P_2}{\varepsilon ^1} - \sum\limits_{j = 1}^m {\left({{\beta _{j1}}{\phi _{j1}} - {\hat\delta_{j1}}} \right)} + \tilde \Omega\\ 
&- {\mu _2}{c_2}{{\tilde x}^{2T}}{{\hat P}_1}{{\tilde x}^2} - {{\tilde \mu }_1}{\varepsilon ^2}^T{{\hat P}_2}{\varepsilon ^2}  - \sum\limits_{i = m + 1}^N {\left( {{\beta _{i3}}{\phi _{i3}} - {\hat \delta_{i3}}} \right)},
\end{aligned}
\end{equation}
where $\tilde \Omega  = {\Omega _1} + {\hat \Omega _1} + \hat M + \tilde M$.

From Eq.(\ref{type2_31}), it has,
\begin{equation}
\begin{aligned}
{{\dot W}_2} \le -{{\mathord{\buildrel{\lower3pt\hbox{$\scriptscriptstyle\smile$}}
\over \beta } }_1}{W_2} + \mathord{\buildrel{\lower3pt\hbox{$\scriptscriptstyle\smile$}}
\over \Omega }
\end{aligned}
\end{equation}
where ${\mathord{\buildrel{\lower3pt\hbox{$\scriptscriptstyle\smile$}}
\over \beta } _1} = \min \{ {\hat \alpha _2}{c_1},{\beta _{j1}},{\tilde \alpha _1},{\mu _2}{c_2},{\beta _{i3}},{\tilde \mu _1}\}$, $\mathord{\buildrel{\lower3pt\hbox{$\scriptscriptstyle\smile$}}
\over \Omega} = \tilde \Omega + \sum\limits_{j = 1}^m {\hat\delta_{j1}} + \sum\limits_{i = {m+1}}^N {\hat\delta _{i3}}$, $j = 1, \cdots ,m$, $i = m + 1, \cdots ,N$.

Then,
$${W_2} \le ({W_2}(0) - \frac{{\mathord{\buildrel{\lower3pt\hbox{$\scriptscriptstyle\smile$}}
\over \Omega} }}{{{{\mathord{\buildrel{\lower3pt\hbox{$\scriptscriptstyle\smile$}}
\over \beta } }_1}}}){e^{ - {{\mathord{\buildrel{\lower3pt\hbox{$\scriptscriptstyle\smile$}}
\over \beta } }_1}t}} + \frac{{\mathord{\buildrel{\lower3pt\hbox{$\scriptscriptstyle\smile$}}
\over \Omega} }}{{{{\mathord{\buildrel{\lower3pt\hbox{$\scriptscriptstyle\smile$}},
\over \beta } }_1}}}$$

Due to $0 < \mathord{\buildrel{\lower3pt\hbox{$\scriptscriptstyle\smile$}}
\over \Omega} ,{\mathord{\buildrel{\lower3pt\hbox{$\scriptscriptstyle\smile$}}
\over \beta } _1}$, when $t \to \infty$, ${e^{ - {{\mathord{\buildrel{\lower3pt\hbox{$\scriptscriptstyle\smile$}}
\over \beta } }_1}t}} \to 0$, that is, for $\forall \varepsilon  > 0$, $\exists T > 0$, when $t > T$, there is ${e^{ - {{\mathord{\buildrel{\lower3pt\hbox{$\scriptscriptstyle\smile$}}
\over \beta } }_1}t}} < \varepsilon$. This shows that ${W_3}$ is bounded.

From Eq.(\ref{type1_12.1}), it has,
\begin{equation}
\begin{array}{l}
||{\tilde x^1}|{|^2} \le \lambda _{\min }^{ - 1}({P_1})V_2^1 \le \lambda _{\min }^{ - 1}({P_1}){W_2},\\
||{\varepsilon ^1}|{|^2} \le \lambda _{\min }^{ - 1}({P_2})V_2^2 \le \lambda _{\min }^{ - 1}({P_2}){W_2}.
\end{array}
\end{equation}

From Eq.(\ref{type2_21.1}), it has,
\begin{equation}
\begin{array}{l}
||{\tilde x^2}|{|^2} \le \lambda _{\min }^{ - 1}({\hat P_1})\hat V_2^1 \le \lambda _{\min }^{ - 1}({\hat P_1}){W_2}, \\
||{\varepsilon ^2}|{|^2} \le \lambda _{\min }^{ - 1}({\hat P_2})\hat V_2^2 \le \lambda _{\min }^{ - 1}({\hat P_2}){W_2}.
\end{array}
\end{equation}

Therefore, $\tilde x_i^1$, $\varepsilon _i^1$, $\tilde x_j^2$ , $\varepsilon _i^2$, $i = 1, \cdots ,m$, $j = m + 1, \cdots ,N$, is UUB. $\hfill{} \Box$

\section{Proof of Lemma 12}

Suppose that the $i^{th}$ agent involves the Zeno-behavior, implying that there exists ${T_0} > 0$ such that $\mathop {\lim }\limits_{k \to \infty } t_k^i = {T_0}$, that is, for $\forall \varepsilon  > 0$, $\exists N(\varepsilon ) > 0$, when $N(\varepsilon ) < k$, it has $t_k^i \in [{T_0} - \varepsilon ,{T_0})$.

From Eq.(\ref{type1_10.2}), it has,
\begin{equation} \label{type2_32}
\begin{array}{l}
{\phi _{i1}} \ge {\phi _{i1}}(0){e^{ - \left( {{\beta _{i1}}{\rm{ + }}\frac{1}{{{\omega _{i1}}}}} \right)t}} > 0.
\end{array}
\end{equation}

From Eq.(\ref{type1_10.1}), it has,
\begin{equation}
\begin{aligned}
||\dot e_{{v_i}}^1|| \le& ||\dot {\bar \zeta} _j^2 - \dot \zeta _j^2|| 
\le ||{A_0}e_{{v_i}}^1|| + ||M_i^{11}M_i^{12}\sum\limits_{j \in \Lambda 1} {L_{ij}^1(\bar \zeta _j^1 - \bar \zeta _i^1)} ||\\
&+ ||M_i^{11}{{\bar g}_{0i}}M_i^1({x_0} - \zeta _i^1)|| + ||M_i^{13}\sum\limits_{j = m + 1}^N {L_{ij}^{12}(\bar \zeta _j^2 - \bar \zeta _i^1)} ||.
\end{aligned}
\end{equation}

According to Eq.(\ref{type1_10.2}) and Lemma 11, $||\dot e_{{v_i}}^1||$ is bounded. Assuming that the boundary of $||\dot e_{{v_i}}^1||$ is ${r_1}$, that is, $||\dot e_{{v_i}}^1|| \le {r_1}$.

Then,
\begin{equation} \label{type2_33}
\begin{array}{l}
||e_{{v_i}}^1 - e_{{v_i}}^1(t_k^i)|| \le {r_1}||t - t_k^i||,
\end{array}
\end{equation}
where $e_{{v_i}}^1(t_k^i) = 0$. Through the triggering mechanism (\ref{type1_10.2}), it has,
\begin{equation}
\begin{array}{l}
{\omega _{i1}}({\gamma _{i1}}(e_{{v_i}}^{1T}{H_{i1}}e_{{v_i}}^1) - \frac{1}{{{\theta _{i1}}}}||{\psi _{i1}}|{|^2} - {\hat\delta_{i1}}) \le {\phi _{i1}}.
\end{array}
\end{equation}

Then,
\begin{equation} \label{type2_34}
\begin{aligned}
||e_{{v_i}}^1|{|^2} \le& \frac{1}{{{\lambda _{\min }}({H_{i1}}){\omega _{i1}}}}{\phi _{i1}}
+ \frac{1}{{{\lambda _{\min }}({H_{i1}}){\gamma _{i1}}{\theta _{i1}}}}||{\psi _{i1}}|{|^2} 
 + \frac{\hat \delta_{i1}}{{{\lambda _{\min }}({H_{i1}}){\gamma _{i1}}}}.
\end{aligned}
\end{equation}

According to Eq.(\ref{type2_32}), Eq.(\ref{type2_34}) is to be established so that the following inequality holds,
\begin{equation} \label{type2_35}
\begin{array}{l}
||e_{{v_i}}^1|| \le {r_2}{e^{ - \frac{1}{2}\left( {{\beta _{i1}}{\rm{ + }}\frac{1}{{{\omega _{i1}}}}} \right)t}},
\end{array}
\end{equation}
where ${r_2}$ is a positive real value.

From Eq.(\ref{type2_33}), the following condition is a sufficient condition to ensure that Eq.(\ref{type2_35}) holds,
\begin{equation} \label{type2_36}
\begin{array}{l}
{r_1}||t - t_k^i|| \le {r_2}{e^{ - \frac{1}{2}\left( {{\beta _{i1}}{\rm{ + }}\frac{1}{{{\omega _{i1}}}}} \right)t}}.
\end{array}
\end{equation}

Define the following triggering function,
\begin{equation} \label{type2_37}
\begin{array}{l}
{\mathord{\buildrel{\lower3pt\hbox{$\scriptscriptstyle\smile$}}
\over h} _i} = ||t - t_k^i|| - \frac{{{r_2}}}{{{r_1}}}{e^{ - \frac{1}{2}\left( {{\beta _{i1}}{\rm{ + }}\frac{1}{{{\omega _{i1}}}}} \right)t}}.
\end{array}
\end{equation}

Considering the following triggering mechanisms,
\begin{equation} \label{type2_38}
\begin{array}{l}
\mathord{\buildrel{\lower3pt\hbox{$\scriptscriptstyle\smile$}}
\over l} _{k + 1}^i = \mathop {\inf }\limits_{\mathord{\buildrel{\lower3pt\hbox{$\scriptscriptstyle\smile$}}
\over r}  \ge \mathord{\buildrel{\lower3pt\hbox{$\scriptscriptstyle\smile$}}
\over t} _k^i} \{ \mathord{\buildrel{\lower3pt\hbox{$\scriptscriptstyle\smile$}}
\over r} :{\mathord{\buildrel{\lower3pt\hbox{$\scriptscriptstyle\smile$}}
\over h} _i} \ge 0\}.
\end{array}
\end{equation}

Let ${\varepsilon _1} = \frac{{2{r_2}}}{{3{r_1}}}{e^{ - \frac{1}{2}\left( {{\beta _{i1}}{\rm{ + }}\frac{1}{{{\omega _{i1}}}}} \right)t}}$.

From Eq.(\ref{type2_38}), when $N(\varepsilon ) < k$, the trigger interval can be set as $\mathord{\buildrel{\lower3pt\hbox{$\scriptscriptstyle\smile$}}
\over t} _{k + 1}^i - \mathord{\buildrel{\lower3pt\hbox{$\scriptscriptstyle\smile$}}
\over t} _k^i = \frac{3}{2}{\varepsilon _1}$.

According to Eqs.(\ref{type2_35}) and (\ref{type2_36}), when $N(\varepsilon ) < k$, the trigger interval of Eq.(\ref{type1_10.2}) is satisfied,
$$t_{k + 1}^i - t_k^i \ge \mathord{\buildrel{\lower3pt\hbox{$\scriptscriptstyle\smile$}}
\over t} _{k + 1}^i - \mathord{\buildrel{\lower3pt\hbox{$\scriptscriptstyle\smile$}}
\over t} _k^i = \frac{3}{2}{\varepsilon _1}.$$

This is in contradiction with $t_k^i \in [{T_0} - \varepsilon ,{T_0})$. Thus, Zeno phenomenon is excluded.       $\hfill{} \Box$

\section*{Acknowledgement} 
This work was supported in part by the National Natural Science Foundation of China Ref. U1913201 and Ref. 61973296, in part by the Guangdong Basic and Applied Basic Research Foundation Ref.2021B1515120038, in part by the Shenzhen Science and Technology Innovation Commission Project Grant Ref.JCYJ20200109114839874, Ref.JSGG20210802154535003 and Ref.JCYJ20210324101215039.

\bibliographystyle{elsarticle-num} 
\bibliography{referen_cy}        % bibliography (preferred). The
                                 % correct style is generated by
                                 % Elsevier at the time of printing.
%\par\noindent
%\parbox[t]{\linewidth}{
%\noindent\parpic{\includegraphics[height=1.2in,width=1in,clip,keepaspectratio]{zheng.jpg}}
%\noindent {\bf Leyi Zheng} received the M.S. degree in basic mathematics, School of Mathematical Sciences from Huaqiao University, China, in 2019. He is currently pursuing the Ph.D. degree with the Shenzhen Institute of Advanced Technology, Chinese Academy of Sciences, Shenzhen, China. His current research interests include Multi-agent cooperative control and Fault tolerant control.}
%
%\vspace{2\baselineskip}
%\vspace{2\baselineskip}
%\par\noindent
%\parbox[t]{\linewidth}{
%\noindent\parpic{\includegraphics[height=1in,width=1in,clip,keepaspectratio]{zhou.png}}
%\noindent {\bf Yimin Zhou} received the Ph.D. degree in control engineering from the University of Oxford, U.K., in 2008. She is currently a Full Professor with the Shenzhen Institute of Advanced Technology, Chinese Academy Sciences, Shenzhen, China. Her research interests include nonlinear control, fault diagnosis, robotics, and energy management.}
%\vspace{2\baselineskip}
                                                                   
\end{document}